\documentclass[nosumlimits,twoside]{amsart}
\usepackage[cp850]{inputenc}
\usepackage[bookmarksnumbered,plainpages,hyperindex]{hyperref}
\input{xy}
 \xyoption{all}
\newtheorem{theorem}{\sc Theorem}[section]
\newtheorem{proposition}[theorem]{\sc Proposition}
\newtheorem{notation}[theorem]{\sc Notation}
\newtheorem{lemma}[theorem]{\sc Lemma}
\newtheorem{corollary}[theorem]{\sc Corollary}
\theoremstyle{definition}
\newtheorem{definition}[theorem]{\sc Definition}
\newtheorem{definitions}[theorem]{\sc Definitions}

\newtheorem{examples}[theorem]{\sc Examples}

\theoremstyle{remark}
\newtheorem{remark}[theorem]{\sc Remark}

\newtheorem{claim}[theorem]{}

\newcommand{\J}[2]{i_{{#1}\wedge_E {#2}}^E}

\def\cmc{{}^C\!\calM^C}

\def\Hu{\mathbf{H}^\bullet}

\def\calM{\mathcal{M}}
\def\ot{\otimes}
\def\cot{\square}

\def\K{\textrm{Ker}}
\def\C{\textrm{Coker}}
\def\Id{\textrm{Id}}
\def\M{\mathcal{M}}

\def\N{\mathbb{N}}
\def\Mt{{^C\mathcal{M}^C}}

\def\w{\wedge_E }

\begin{document}
\pagestyle{headings}
\title{Wedge Products and Cotensor Coalgebras in Monoidal Categories}
\author{A. Ardizzoni}
\address{University of Ferrara, Department of Mathematics, Via Machiavelli
35, Ferrara, I-44100, Italy} \email{alessandro.ardizzoni@unife.it}
\urladdr{http://www.unife.it/utenti/alessandro.ardizzoni}
\thanks{This paper was written while A. Ardizzoni was member of G.N.S.A.G.A. with partial financial support
from M.I.U.R..}
\begin{abstract}
The construction of the cotensor coalgebra for an "abelian
monoidal" category $\M$ which is also cocomplete, complete and
AB5, was performed in [A. Ardizzoni, C. Menini and D. \c{S}tefan,
\emph{Cotensor Coalgebras in Monoidal Categories}, Comm. Algebra,
to appear]. It was also proved that this coalgebra satisfies a
meaningful universal property which resembles the classical one.
Here the lack of the coradical filtration for a coalgebra $E$ in
$\M$ is filled by considering a direct limit $\widetilde{D}$ of a
filtration consisting of wedge products of a subcoalgebra $D$ of
$E$. The main aim of this paper is to characterize hereditary
coalgebras $\widetilde{D}$, where $D$ is a coseparable coalgebra
in $\M$, by means of a cotensor coalgebra: more precisely, we
prove that, under suitable assumptions, $\widetilde{D}$ is
hereditary if and only if it is formally smooth if and only if it
is the cotensor coalgebra $T^c_{D}(D\w D/D)$ if and only if it is
a cotensor coalgebra $T^c_{D}(N)$, where $N$ is a certain
$D$-bicomodule in $\M$. Because of our choice, even when we apply
our results in the category of vector spaces, new results are
obtained.
\end{abstract}
 \maketitle
\section*{Introduction} \markboth{\sc{A. ARDIZZONI}} {\sc{Wedge Products and Cotensor Coalgebras in
Monoidal Categories}}

Let $C$ be a coalgebra over a field $K$ and let $M$ be a
$C$-bicomodule. The cotensor coalgebra $ T^c_C(M)$ was introduced
by Nichols in \cite{Ni} and it appears as a main step in the
classification of finite dimensional Hopf algebras problem (see,
e.g., \cite{AG} and \cite{AS}). In \cite{Ch} the relation between
quiver coalgebras and hereditary coalgebras is investigated. In
\cite{JLMS}, hereditary coalgebras with coseparable coradical are
characterized by means of a suitable cotensor coalgebra: more
precisely, the authors prove that a coalgebra $C$ with coseparable
coradical $D$ is hereditary if and only if it is formally smooth
if and only if it is a cotensor coalgebra $T^c_{D}(N)$, where $N$
is a certain $D$-bicomodule.
\newline The main aim of this paper is to prove Theorem \ref{teo: hered <=> fs} which
establishes an analogous result inside the framework of monoidal
categories. This is pursued using the notion of formally smooth
coalgebra for "abelian monoidal" categories developed in
\cite{AMS} and using the construction of the cotensor coalgebra
for "abelian monoidal" categories  performed in \cite{Cotensor}
.\newline The basic point when dealing with coalgebras in monoidal
categories is that there is no notion of coradical. The idea then
is to take a subcoalgebra $D$ of a coalgebra $C$ and to consider
the coalgebra $\widetilde{D}$ which is the direct limit of the
iterated wedge powers of $D$ in $E$. Then the coalgebra $D$ acts,
in a certain sense, as the coradical of $\widetilde{D}$. Thus,
because of our choice, even when we apply our results (e.g.
Theorem \ref{teo: Teo2}) in the category of vector spaces, new
results are obtained. It is also interesting to point out that,
working in this wider context, we had to develop some properties
of the wedge product that have an intrinsic interest. Due to the
width of our setting, many technical results were needed. To
enable an easier reading, we decided to postpone a number of them
in two appendices that can be found at the end of the
paper.\newline \medskip

\textbf{Notations.} \ Let $[(X,i_{X})]$ be a subobject of an
object $E$ in an abelian category $\M,$ where
$i_{X}=i_{X}^{E}:X\hookrightarrow E$ is a monomorphism and
$[(X,i_{X})]$ is the associated equivalence class. By abuse of
language, we will say that $(X,i_{X})$ is a subobject of $E$ and
we will write $(X,i_{X})=(Y,i_{Y})$ to mean that $(Y,i_{Y}) \in
[(X,i_{X})]$. The same convention applies to cokernels. If
$(X,i_{X})$ is a subobject of $E$ then we will write
$(E/X,p_X)=\C(i_X)$, where $p_{X}=p_{X}^{E}:E\rightarrow
E/X$.\newline
Let $(X_{1},i_{X_{1}}^{Y_{1}})$ be a subobject of $Y_{1}$ and let $%
(X_{2},i_{X_{2}}^{Y_{2}})$ be a subobject of $Y_{2}$. Let $%
x:X_{1}\rightarrow X_{2}$  and $y:Y_{1}\rightarrow Y_{2}$ be
morphisms such that $y\circ
i_{X_{1}}^{Y_{1}}=i_{X_{2}}^{Y_{2}}\circ x$. Then there exists
a unique morphism, which we denote by $y/x=\frac{y}{x}:Y_{1}/{X_{1}}%
\rightarrow Y_{2}/{X_{2},}$ such that $\frac{y}{x}\circ
p_{X_{1}}^{Y_{1}}=p_{X_{2}}^{Y_{2}}\circ y$:
\begin{equation*}
\xymatrix@R=20pt@C=40pt{
  X_1  \ar[d]_{x} \ar@{^{(}->}[r]^{i_{X_1}^{Y_1}} & Y_1 \ar[d]_{y} \ar[r]^{p_{X_1}^{Y_1}} & \frac{Y_1}{X_1} \ar[d]^{\frac{y}{x}} \\
  X_2  \ar@{^{(}->}[r]^{i_{X_2}^{Y_2}} & Y_2 \ar[r]^{p_{X_2}^{Y_2}} & \frac{Y_2}{X_2}   }
\end{equation*}
\section{Monoidal Categories}

\begin{claim}
\label{MonCat}A \emph{monoidal category} means a category
$\mathcal{M}$ that is endowed with a functor $\otimes
:\mathcal{M}\times \mathcal{M}\rightarrow \mathcal{M}$, an object
$\mathbf{1}\in \mathcal{M}$ and functorial isomorphisms:
$a_{X,Y,Z}:(X\otimes Y)\otimes Z\rightarrow X\otimes (Y\otimes
Z),$ $l_{X}:\mathbf{1}\otimes X\rightarrow X$ and $r_{X}:X\otimes \mathbf{1}%
\rightarrow X.$ The functorial morphism $a$ is called the \emph{%
associativity constraint }and\emph{\ }satisfies the \emph{Pentagon Axiom, }%
that is the following diagram
\begin{equation*}
                  \xymatrix@R=45pt@C=-30pt{
                  &((U\otimes V)\otimes W)\otimes X
                  \ar[rr]^{ \alpha_{U,V,W}\otimes X}
                  \ar[dl]|{ \alpha_{U\otimes V,W,X}}
                  &&(U\otimes (V\otimes W))\otimes X
                  \ar[dr]|{ \alpha_{U,V\otimes W,X}}&\\
                  (U\otimes V)\otimes (W\otimes X)
                  \ar[drr]|{ \alpha_{U,V,W\otimes X}}
                  &&&&U\otimes ((V\otimes W)\otimes X)
                  \ar[dll]|{ U\otimes \alpha_{V,W,X}}
                  \\&&U\otimes (V\otimes (W\otimes X))&&&
                 }
\end{equation*}
is commutative, for every $U,\,V,$ $W,$ $X$ in $\mathcal{M}.$ The morphisms $%
l$ and $r$ are called the \emph{unit constraints} and they are
assumed to satisfy the \emph{Triangle Axiom, }i.e. the following
diagram
\begin{equation*}\xymatrix@R=30pt@C=-2pt{
                 (V\ot \mathbf{1})\ot W \ar[dr]_{r_V\ot
                 W}\ar[rr]^{a_{V,\mathbf{1},W}}&&V\ot (\mathbf{1}\ot
                 W)\ar[dl]^{V\ot l_W} &\\
                 &V\ot W
                 }
\end{equation*}
is commutative. The object $\mathbf{1}$ is called the \emph{unit} of $%
\mathcal{M}$. For details on monoidal categories we refer to \cite[Chapter XI]%
{Ka} and \cite{Maj2}. A monoidal category is called \emph{strict}
if the associativity constraint and unit constraints are the
corresponding identity morphisms.
\end{claim}

\begin{claim}
\label{cl:CohThm}As it is noticed in \cite[p. 420]{Maj2}, the
Pentagon Axiom solves the consistency problem that appears because
there are two ways to go from $((U\otimes V)\otimes W)\otimes X$
to $U\otimes (V\otimes (W\otimes X)).$ The coherence theorem, due
to S. Mac Lane, solves the similar problem for the tensor product
of an arbitrary number of objects in $\mathcal{M}.$ Accordingly
with this theorem, we can always omit all brackets and simply
write $X_{1}\otimes \cdots \otimes X_{n}$ for any object obtained from $%
X_{1},\ldots ,X_{n}$ by using $\otimes $ and brackets. Also as a
consequence of the coherence theorem, the morphisms $a,$ $l,$ $r$
take care of themselves, so they can be omitted in any computation
involving morphisms in $\mathcal{M.}$\newline The notions of
algebra, module over an algebra, coalgebra and comodule over a
coalgebra can be introduced in the general setting of monoidal
categories. For more details, see \cite{AMS}.
\end{claim}

\section{Wedge and cotensor products}

We quote from \cite[2.4]{Cotensor}

\begin{definition}\label{abelianmonoidal}A monoidal category
$(\M,\ot,\mathbf{1})$ will be called an \emph{\textbf{abelian
monoidal category} } if:
\begin{enumerate}
    \item $\M$ is an abelian category
    \item both the functors $X\ot (-):\M\to\M$ and $(-)\ot
    X:\M\to\M$ are additive and left exact, for every object $X\in
    \M$.
\end{enumerate}
\end{definition}

\begin{claim}
Let $E$ be a coalgebra in an abelian monoidal category
$\mathcal{M}$. Let us recall,
(see \cite[page 60]{Mo}), the definition of wedge of two subobjects $%
X,Y $ of $E$ in $\M:$%
\begin{equation*}
(X\wedge_E Y,i_{X\w Y}^E):=Ker[ (p _{X}\otimes p _{Y}) \circ
\triangle _{E}] ,
\end{equation*}
where $p _{X}:E\rightarrow E/X$ and $p _{Y}:E\rightarrow E/Y$ are
the canonical quotient maps. In particular we have the following
exact sequence:
\begin{equation*}
\xymatrix@C=0.9cm{
  0 \ar[r] & X\wedge_E Y \ar[rr]^{i_{X\w Y}^E} && E \ar[rr]^(.4){(p _{X}\otimes p _{Y}) \circ \triangle
  _{E}} && E/X \ot E/Y.}
\end{equation*}
Consider the following commutative diagrams in $\M$
  \begin{equation*}
\begin{tabular}{cc}
$\xymatrix@C=2cm{
  X_1 \ar[d]_{x} \ar@{^{(}->}[r]^{i_{X_1}^{E_1}}
                & E_1 \ar[d]^{e}  \\
  X_2 \ar@{^{(}->}[r]_{i_{X_2}^{E_2}}
                & E_2             }
$& $ \xymatrix@C=2cm{
  Y_1 \ar[d]_{y} \ar@{^{(}->}[r]^{i_{Y_1}^{E_1}}
                & E_1 \ar[d]^{e}  \\
  Y_2 \ar@{^{(}->}[r]_{i_{Y_2}^{E_2}}
                & E_2             }
$
\end{tabular}
\end{equation*}  where $e$ is a coalgebra homomorphism. Then there is a unique
morphism $x\wedge_e y:X_1\wedge_{E_1} Y_1\to X_2\wedge_{E_2} Y_2$
such that the following diagram
\begin{equation*}
\xymatrix@C=2cm{
  X_1\wedge_{E_1} Y_1 \ar@{.>}[d]_{x\wedge_e y} \ar[r]^{i_{X_1\wedge_{E_1} Y_1}^{E_1}}
                & E_1 \ar[d]^{e}  \\
  X_2\wedge_{E_2} Y_2 \ar[r]_{i_{X_2\wedge_{E_2} Y_2}^{E_2}}
                & E_2             }
\end{equation*}
commutes. In fact we have
\begin{eqnarray*}
&&(p_{X_{2}}^{E_{2}}\otimes p_{Y_{2}}^{E_{2}})\circ \Delta
_{E_{2}}\circ
e\circ i_{X_1\wedge_{E_1} Y_1}^{E_1} \\
&=&(p_{X_{2}}^{E_{2}}\otimes p_{Y_{2}}^{E_{2}})\circ (e\otimes
e)\circ \Delta
_{E_{1}}\circ i_{X_1\wedge_{E_1} Y_1}^{E_1} \\
&=&(\frac{e}{x}\otimes \frac{e}{y})\circ (p_{X_{1}}^{E_{1}}\otimes
p_{Y_{1}}^{E_{1}})\circ \Delta _{E_{1}}\circ i_{X_1\wedge_{E_1}
Y_1}^{E_1}=0
\end{eqnarray*}
so that, since $(X_2\wedge_{E_2} Y_2,i_{X_2\wedge_{E_2}
Y_2}^{E_2})$ is the kernel of $(p_{X_{2}}^{E_{2}}\otimes
p_{Y_{2}}^{E_{2}})\circ \Delta _{E_{2}}$, we conclude.
\end{claim}

\begin{lemma}
Consider the following commutative diagrams in $\M$
  \begin{equation*}
\begin{tabular}{cc}
$\xymatrix@C=2cm{
  X_1 \ar[d]_{x} \ar@{^{(}->}[r]^{i_{X_1}^{E_1}}
                & E_1 \ar[d]^{e}  \\
  X_2 \ar[d]_{x'}\ar@{^{(}->}[r]_{i_{X_2}^{E_2}}
                & E_2 \ar[d]_{e'}           \\
  X_3 \ar@{^{(}->}[r]_{i_{X_3}^{E_3}}
                & E_3 }
$& $\xymatrix@C=2cm{
  Y_1 \ar[d]_{y} \ar@{^{(}->}[r]^{i_{Y_1}^{E_1}}
                & E_1 \ar[d]^{e}  \\
  Y_2 \ar[d]_{y'}\ar@{^{(}->}[r]_{i_{Y_2}^{E_2}}
                & E_2 \ar[d]_{e'}           \\
  Y_3 \ar@{^{(}->}[r]_{i_{Y_3}^{E_3}}
                & E_3 }$
\end{tabular}
\end{equation*}  where $e$ and $e'$ are coalgebra
homomorphisms. Then we have
\begin{equation}\label{formula compos of wedge}
  (x'\wedge_{e'} y')\circ  (x\wedge_{e} y)= (x'x\wedge_{e'e} y'y)
\end{equation}
\end{lemma}

\begin{proof}: straightforward.
\end{proof}

\begin{claim}\label{claim: cotensor} Let $\mathcal{M}$ be an abelian monoidal category and let $E$ be a
coalgebra in $(\M,\ot,\mathbf{1})$. Given a right $E$-bicomodule
$(V,\rho^r_V)$ and a left $E$-comodule $(W,\rho^l_W)$, their
cotensor product over $E$ in $\mathcal{M}$ is defined to be the
equalizer $(V\cot_{E}W,{\chi(V,W)=\chi_E(V,W)})$ of the couple of
morphism $(\rho^r_V\ot W,V\ot \rho^l_W)$:
\begin{equation*}
\xymatrix@C=1cm{
  0 \ar[r] & V\cot_{E}W \ar[rr]^{\chi_E(V,W)} && V\otimes W \ar@<.5ex>[rr]^{\rho^r_V\ot W} \ar@<-.5ex>[rr]_{V\ot \rho^l_W}&&V\ot E\ot W  }
\end{equation*}
Since the tensor functors are left exact, in view of
\cite[Proposition 1.3]{Cotensor}, then $V\Box _{E}W$ is also a
$E$-bicomodule, namely it is $E$-sub-bicomodule of $V\ot W$,
whenever $V$ and $W$ are $E$-bicomodules. Furthermore, in this
case, the category $({^E\M^E},\cot_E,E)$ is still an abelian
monoidal category; the associative and unit constraints are
induced by the ones in $\M$ (the proof is dual to \cite[Theorem
1.11]{AMS}). Therefore, also using $\cot_E$, one can forget about
brackets. Moreover the functors $M\cot_E(-):{^E\M}\to \M$ and
$(-)\cot_E M:{\M^E}\to\M$ are left exact for any $M\in\M.$

We will write $\cot$ instead of $\cot_C$, whenever there is no
danger of misunderstanding.\\One has the
following result.\end{claim}

\begin{claim}Let $e:E_{1}\rightarrow E_{2}$ be a coalgebra homomorphism in an abelian monoidal category $\M$. Let $%
(V_{1},\rho _{V_{1}}^{E_{1}})$ be a right $E_{1}$-comodule, let $(W_{1},{%
^{E_{1}}\rho _{W_{1}}})$ be a left $E_{1}$-comodule, let
$(V_{2},\rho _{V_{2}}^{E_{2}})$ be a right $E_{2}$-comodule and
let $(W_{2},{^{E_{2}}\rho
_{W_{2}}})$ be a left $E_{2}$-comodule. Let $v:V_{1}\rightarrow V_{2}$ and $%
w:W_{1}\rightarrow W_{2}$ be $E_{2}$-comodule homomorphisms (where
$V_{1}$ and $W_{1}$ are regarded as $E_{2}$-comodules via $e$).
Then there is a unique morphism $v\cot _{e}w:V_{1}\cot
_{E_{1}}W_{1}\rightarrow V_{2}\cot _{E_{2}}W_{2}$ such that the
following diagram
\begin{equation*}
  \xymatrix@C=2cm{
  V_1\cot_{E_1} W_1 \ar@{.>}[d]_{v\cot_e w} \ar[r]^{\chi_{E_1}(V_1,W_1)}
                & V_1\ot W_1 \ar[d]^{v\ot w}  \\
  V_2\cot_{E_2} W_2 \ar[r]_{\chi_{E_2}(V_2,W_2)}
                & V_2\ot W_2             }
\end{equation*}
commutes. In fact we have
\begin{eqnarray*}
&&(\rho _{V_{2}}^{E_{2}}\otimes W_{2})\circ (v\otimes w)\circ \chi
_{E_{1}}(V_{1},W_{1}) \\
&=&(v\otimes E_{2}\otimes w)\circ (\rho _{V_{1}}^{E_{2}}\otimes
W_{1})\circ
\chi _{E_{1}}(V_{1},W_{1}) \\
&=&(v\otimes E_{2}\otimes w)\circ \lbrack (V_{1}\otimes e)\rho
_{V_{1}}^{E_{1}}\otimes W_{1}]\circ \chi _{E_{1}}(V_{1},W_{1}) \\
&=&(v\otimes e\otimes w)\circ (\rho _{V_{1}}^{E_{1}}\otimes
W_{1})\circ \chi
_{E_{1}}(V_{1},W_{1}) \\
&=&(v\otimes e\otimes w)\circ \lbrack V_{1}\otimes {^{E_{1}}\rho _{W_{1}}}%
]\circ \chi _{E_{1}}(V_{1},W_{1}) \\
&=&(v\otimes E_{2}\otimes w)\circ \lbrack V_{1}\otimes (e\otimes
W_{1})\circ
{^{E_{1}}\rho _{W_{1}}]}\circ \chi _{E_{1}}(V_{1},W_{1}) \\
&=&(v\otimes E_{2}\otimes w)\circ (V_{1}\otimes {^{E_{2}}\rho _{W_{1}})}%
\circ \chi _{E_{1}}(V_{1},W_{1}) \\
&=&(V_{2}\otimes {^{E_{2}}\rho _{W_{2}}})\circ (v\otimes w)\circ
\chi _{E_{1}}(V_{1},W_{1})
\end{eqnarray*}
so that, since $(V_{2}\cot _{E_{2}}W_{2},\chi
_{E_{2}}(V_{2},W_{2}))$ is the
equalizer of $\rho _{V_{2}}^{E_{2}}\otimes W_{2}$ and $V_{2}\otimes {%
^{E_{2}}\rho _{W_{2}},}$ we conclude.\\
Note that if $E_1=E_2=E$ and $e=\Id_E$, one has $$v\cot_e
w=v\cot_E w.$$
\end{claim}

\begin{lemma}
Let $e:E_{1}\rightarrow E_{2}$ and $e^{\prime }:E_{2}\rightarrow
E_{3}$ be coalgebra homomorphisms in $\M$. Let
\begin{eqnarray*}
(V_{1},\rho _{V_{1}}^{E_{1}}) \in {\M^{E_{1}},}\text{\quad
}(V_{2},\rho _{V_{2}}^{E_{2}})\in {\M^{E_{2}},}\text{\quad
}(V_{3},\rho
_{V_{3}}^{E_{3}})\in {\M^{E_{3}},} \\
(W_{1},{^{E_{1}}\rho _{W_{1}}}) \in {^{E_{1}}\M,}\text{\quad }%
(W_{2},^{E_{2}}\rho _{W_{2}})\in {^{E_{2}}\M,}\text{\quad }%
(W_{3},^{E_{3}}\rho _{W_{3}})\in {^{E_{3}}\M.}
\end{eqnarray*}
Let $v:V_{1}\rightarrow V_{2}$ and $w:W_{1}\rightarrow W_{2}$ be $E_{2}$%
-comodule homomorphisms (where $V_{1}$ and $W_{1}$ are regarded as $E_{2}$%
-comodules via $e$) and let $v^{\prime }:V_{2}\rightarrow V_{3}$ and $%
w^{\prime }:W_{2}\rightarrow W_{3}$ be $E_{3}$-comodule
homomorphisms (where $V_{2}$ and $W_{2}$ are regarded as
$E_{3}$-comodules via $e^{\prime }$). Then
\begin{equation}\label{formula compos of cotensor}
(v^{\prime }\cot_{e^{\prime }}w^{\prime })\circ (v\cot_{e}w)=(v^{\prime }v%
\cot_{e^{\prime }e}w^{\prime }w).
\end{equation}
\end{lemma}

\begin{proof}
: straightforward.
\end{proof}

\begin{proposition}\label{pro: rho bar is iso}Let $\mathcal{M}$ be an abelian monoidal
category, let $E$ be a coalgebra in $\M$ and let $(V,\rho^E_V)$ be
a right $E$-comodule. Then the morphism $\rho
_{V}^{E}:V\rightarrow V\ot E$ factorizes to a unique
morphism $\overline{\rho }_{V}^{E}:V\rightarrow V\square _{E}E$ such that $%
\chi (V,E)\circ \overline{\rho }_{V}^{E}=\rho _{V}^{E}.$ Moreover $\overline{%
\rho }_{V}^{E}$ is an isomorphism whose inverse is given by
\begin{equation*}
(\overline{\rho }_{V}^{E})^{-1}=r_{V}\circ (V\ot \varepsilon
_{E})\circ \chi (V,E).
\end{equation*}
An analogous statement holds for a left $E$-comodule $(W,\rho
_{W}^{l})$.
\end{proposition}

\begin{proof}
 Since $(V,\rho^E_V)$ is a right $E$-comodule, by universal property of the equalizer, we get the existence of
 $\overline{\rho }_{V}^{E}$.  %
 Set $s_{V}:=r_{V}\circ (V\ot \varepsilon
_{E})\circ \chi (V,E).$ We have
\begin{equation*}
s_{V}\circ \overline{\rho }_{V}^{E}=r_{V}\circ (V\ot
\varepsilon _{E})\circ \chi (V,E)\circ \overline{\rho }_{V}^{E}=r_{V}\circ (V%
\ot \varepsilon _{E})\circ \rho _{V}^{E}=r_{V}\circ
r_{V}^{-1}=\mathrm{Id}_{V}.
\end{equation*}
Moreover we have
\begin{eqnarray*}
\chi (V,E)\circ \overline{\rho }_{V}^{E}\circ s_{V} &=&\rho
_{V}^{E}\circ
s_{V} \\
&=&\rho _{V}^{E}\circ r_{V}\circ (V\otimes \varepsilon _{E})\circ
\chi (V,E)
\\
&=&r_{V\otimes E}\circ (\rho _{V}^{E}\otimes \mathbf{1})\circ
(V\otimes
\varepsilon _{E})\circ \chi (V,E) \\
&=&r_{V\otimes E}\circ (V\otimes E\otimes \varepsilon _{E})\circ
(\rho
_{V}^{E}\otimes E)\circ \chi (V,E) \\
&=&r_{V\otimes E}\circ (V\otimes E\otimes \varepsilon _{E})\circ
(V\otimes
\Delta _{E})\circ \chi (V,E) \\
&=&(V\otimes r_{E})\circ \lbrack V\otimes (E\otimes \varepsilon
_{E})\circ
\Delta _{E}]\circ \chi (V,E) \\
&=&(V\otimes r_{E})\circ (V\otimes r_{E}^{-1})\circ \chi
(V,E)=\chi (V,E).
\end{eqnarray*}
Since $\chi (V,E)$ is a monomorphism, we get $\overline{\rho
}_{V}^{E}\circ s_{V}=\mathrm{Id}_{V\otimes E}.$
\end{proof}

\begin{lemma}
\label{lem: induced comodule} Let $\alpha :A\rightarrow E$ be a
monomorphism
which is a coalgebra homomorphism in an abelian monoidal category $\mathcal{M%
}$. Let $(W,{^{E}\rho _{W}})$ be a left $E$%
-comodule and let $^{A}\rho _{W}:W\rightarrow A\otimes W$ be a
morphism such that
\begin{equation*}
^{E}\rho _{W}=(\alpha \otimes W)\circ {^{A}\rho _{W}.}
\end{equation*}
Then $(W,{^{A}\rho _{W}})$ is a left $A$-comodule. \newline
Let ${^{A}}\overline{{\rho }}{_{W}}$ be the unique morphism such that ${%
^{A}\rho _{W}=\chi (A,W)\circ {^{A}}\overline{{\rho }}{_{W}}}$. Then ${^{A}}%
\overline{{\rho }}{_{W}:W\rightarrow A}\cot_{A}W$ is a morphism of left $E$%
-comodules.
\end{lemma}

\begin{proof}
We have
\begin{eqnarray*}
&&(\alpha \otimes \alpha \otimes W)\circ (A\otimes {^{A}\rho _{W}})\circ {%
^{A}\rho _{W}} \\
&=&[E\otimes (\alpha \otimes W){^{A}\rho _{W}}]\circ (\alpha
\otimes W)\circ
{^{A}\rho _{W}} \\
&=&\left( E\otimes {^{E}\rho _{W}}\right) \circ {^{E}\rho _{W}} \\
&=&\left( \Delta _{E}\otimes W\right) \circ {^{E}\rho _{W}} \\
&=&\left( \Delta _{E}\otimes W\right) \circ (\alpha \otimes W)\circ {%
^{A}\rho _{W}}=(\alpha \otimes \alpha \otimes W)\circ \left(
\Delta _{A}\otimes W\right) \circ {^{A}\rho _{W}.}
\end{eqnarray*}
Since $\alpha \otimes \alpha \otimes W$ is a monomorphism, we get
\begin{equation*}
(A\otimes {^{A}\rho _{W}})\circ {^{A}\rho _{W}=}\left( \Delta
_{A}\otimes W\right) \circ {^{A}\rho _{W}.}
\end{equation*}
Moreover we have
\begin{equation*}
(\varepsilon _{A}\otimes W)\circ {^{A}\rho _{W}}=(\varepsilon
_{E}\otimes W)\circ (\alpha \otimes W)\circ {^{A}\rho
_{W}}=(\varepsilon _{E}\otimes W)\circ {^{E}\rho _{W}}=l_{W}.
\end{equation*}
Therefore $(W,{^{A}\rho _{W}})$ is a left $A$-comodule. Let us prove that ${^{A}}\overline{{\rho }}{_{W}}$ is a morphism of left $E$%
-comodules. We have:
\begin{eqnarray*}
&&(E\otimes \alpha \otimes W)\circ \lbrack E\otimes \chi
(A,W)]\circ
(E\otimes {^{A}}\overline{{\rho }}{_{W}})\circ {^{E}\rho _{W}} \\
&=&(E\otimes \alpha \otimes W)\circ (E\otimes {^{A}\rho _{W}})\circ {%
^{E}\rho _{W}} \\
&=&(E\otimes {^{E}\rho _{W}})\circ {^{E}\rho _{W}} \\
&=&\left( \Delta _{E}\otimes W\right) \circ {^{E}\rho _{W}} \\
&=&\left( \Delta _{E}\otimes W\right) \circ (\alpha \otimes W)\circ {%
^{A}\rho _{W}} \\
&=&\left( \Delta _{E}\otimes W\right) \circ (\alpha \otimes
W)\circ \chi
(A,W)\circ {^{A}\overline{{\rho }}_{W}} \\
&=&(E\otimes \alpha \otimes W)\circ \left( {^{E}\rho _{A}}\otimes
W\right)
\circ \chi (A,W)\circ ^{A}\overline{{\rho }}_{W} \\
&=&(E\otimes \alpha \otimes W)\circ \lbrack E\otimes \chi
(A,W)]\circ \left(
{^{E}\rho _{A}}\cot_{A}W\right) \circ {^{A}\overline{{\rho }}_{W}} \\
&=&(E\otimes \alpha \otimes W)\circ \lbrack E\otimes \chi (A,W)]\circ {%
^{E}\rho _{A\cot_{A}W}}\circ {^{A}\overline{{\rho }}_{W}}.
\end{eqnarray*}
Since $E\otimes \alpha \otimes W$ and $E\otimes \chi (A,W)$ are
monomorphisms, we obtain:
\begin{equation*}
(E\otimes {^{A}}\overline{{\rho }}{_{W}})\circ {^{E}\rho _{W}={^{E}\rho _{A%
\cot_{A}W}}}\circ ^{A}\overline{{\rho }}_{W},
\end{equation*}
i.e. that ${^{A}}\overline{{\rho }}{_{W}}$ is a morphism of left $E$%
-comodules.
\end{proof}

\begin{remark}
\label{rem utile 1} Let $E$ be a coalgebra in an abelian monoidal category $%
\mathcal{M}$, let $X$ be a right coideal and let $Y$ be a left
coideal of $E$ in $\mathcal{M}$. Then we have
\begin{eqnarray*}
(X\wedge _{E}Y,i_{X\wedge _{E}Y}^{E}) &=&\text{Ker}[(p_{X}\otimes
p_{Y})\circ \triangle _{E}] \\
&=&\text{Ker}[(p_{X}\otimes p_{Y})\circ \chi _{E}(E,E)\circ \overline{%
\triangle }_{E}] \\
&=&\text{Ker}[\chi _{E}(E/X,E/Y)\circ (p_{X}\square
_{E}p_{Y})\circ \overline{\triangle
}_{E}]=\text{Ker}[(p_{X}\square _{E}p_{Y})\circ
\overline{\triangle }_{E}],
\end{eqnarray*}%
where $\overline{\triangle }_{E}:E\rightarrow E\square _{E}E$
denotes the canonical isomorphism. In particular we have the
following exact sequence:
\begin{equation}\label{sequence X wedge E}
\xymatrix@C=0.9cm{
  0 \ar[r] & X\wedge_E Y \ar[rr]^{i_{X\w Y}^E} && E \ar[rr]^(.4){(p _{X}\cot_E p _{Y}) \circ
  \overline{\triangle}
  _{E}} && E/X \cot_E E/Y.}
\end{equation}
\end{remark}

\begin{definition}
Let $\M$ be an abelian monoidal category. %
Let $i_F=i_F^E:F\hookrightarrow E$, $i_A^B:A\hookrightarrow B$ and $i_B=i_B^E:B\hookrightarrow E$ be monomorphisms which are coalgebra homomorphisms.\\
Consider in $\M$ the following commutative diagram with exact rows
and columns
\begin{equation*}
\xymatrix@R=35pt@C=70pt{
                          && F\w A \ar[d]_{i_{F\w A}^E} \ar[r]^{i_{F\w A}^{F\w B}} & F\w B \ar[d]_{i_{F\w B}^E} \\
  &0 \ar[d]_{0} \ar[r]^{0} & E \ar[d]|{(p_{F}\cot _{E}p_{A})\circ \overline{\triangle }_{E}} \ar[r]^{\Id} & E \ar[d]|{(p_{F}\cot _{E}p_{B})\circ \overline{\triangle }_{E}}\ar[r]&0\\
  0\ar[r]&\frac{E}{F}\cot_E \frac{B}{A} \ar[d]_{\Id} \ar[r]^{\frac{E}{F} \cot_E i_{B/A}} & \frac{E}{F} \cot_E \frac{E}{A} \ar[r]^{\frac{E}{F} \cot_E p_{B/A}} & \frac{E}{F} \cot_E \frac{E}{B}\\
  &\frac{E}{F} \cot_E \frac{B}{A}  }
\end{equation*}
In this particular case, we denote the connecting homomorphism by
$$\eta(F,B,A)=\eta^E(F,B,A):F\w B\to \frac{E}{F} \cot_E \frac{B}{A}.$$ Note that, by Proposition \ref{pro: A_1=0}, the morphism $\eta(F,B,A)$ is uniquely
defined by the following relation:
\begin{equation}\label{formula eta}
  (\frac{E}{F} \cot_E i_{B/A}) \circ \eta(F,B,A)=(p_{F}\cot _{E}p_{A})\circ \overline{\triangle }_{E}\circ i_{F\w B}^E.
\end{equation}
By the Snake Lemma we get the following exact sequence:
\begin{equation}\label{sequence eta}
\xymatrix@C=0.9cm{
  0 \ar[r] & F\w A \ar[rr]^{i_{F\w A}^{F\w B}} && F\w B\ar[rr]^{\eta(F,B,A)}&&\frac{E}{F}\cot_E \frac{B}{A} \ar[r]&\C[{(p _{F}\cot_E p _{A}) \circ
  \overline{\triangle}_{E}}].}
\end{equation}
Observe that, if $A=0,$ then $p_A=\Id_E$, so that
$$(p_{F}\cot _{E}p_{A})\circ \overline{\triangle}_{E}=(p_{F}\cot
_{E}E)\circ \overline{\triangle}_{E}=\overline{\rho}^r_{E/F}\circ
p_F$$ is an epimorphism and $$F\w A =\K[(p_{F}\cot _{E}p_{A})\circ
\overline{\triangle}_{E}]=\K(\overline{\rho}^r_{E/F}\circ
p_F)=\K(p_F)=F,$$ that is we have the following exact sequence:
\begin{equation}\label{sequence eta zero}
\xymatrix@C=0.9cm{
  0 \ar[r] & F \ar[rr]^{i_{F}^{F\w B}} && F\w B\ar[rr]^{\eta(F,B,0)}&&\frac{E}{F}\cot_E B \ar[r]&0.}
\end{equation}
\end{definition}

\begin{proposition}
Let $\M$ be an abelian monoidal category. Let $%
i_{F}=i_{F}^{E}:F\hookrightarrow E$, $i_{A}^{B}:A\hookrightarrow B$ and $%
i_{B}=i_{B}^{E}:B\hookrightarrow E$ be monomorphisms which are
coalgebra homomorphisms.\newline Then we have
\begin{equation}\label{formula eta ridotta}
(\frac{E}{F}\cot _{E}\frac{i_{B}^{F\wedge _{E}B}}{A})\circ \eta
(F,B,A)=(p_{F}\circ \J{F}{B}\cot _{E}p_{A}^{F\wedge _{E}B})\circ \overline{%
\triangle }_{F\wedge _{E}B}.
\end{equation}
\end{proposition}

\begin{proof}
By (\ref{formula eta}) we have
\begin{equation*}
(\frac{E}{F}\cot _{E}\frac{i_{B}^{E}}{A})\circ \eta
(F,B,A)=(p_{F}^{E}\cot _{E}p_{A}^{E})\circ \overline{\triangle
}_{E}\circ \J{F}{B}.
\end{equation*}
Thus we obtain
\begin{eqnarray*}
&&(\frac{E}{F}\cot _{E}\frac{i_{F\wedge _{E}B}^{E}}{A})\circ (\frac{E}{F}%
\cot _{E}\frac{i_{B}^{F\wedge _{E}B}}{A})\circ \eta (F,B,A) \\
&=&(\frac{E}{F}\cot _{E}\frac{i_{F\wedge _{E}B}^{E}}{A}\circ \frac{%
i_{B}^{F\wedge _{E}B}}{A})\circ \eta (F,B,A) \\
&=&(\frac{E}{F}\cot _{E}\frac{i_{B}^{E}}{A})\circ \eta (F,B,A) \\
&=&(p_{F}^{E}\cot _{E}p_{A}^{E})\circ \overline{\triangle
}_{E}\circ \J{F}{B}
\\
&=&(p_{F}^{E}\circ \J{F}{B}\cot _{E}p_{A}^{E}\circ \J{F}{B})\circ \overline{%
\triangle }_{F\wedge _{E}B} \\
&=&(p_{F}^{E}\circ \J{F}{B}\cot _{E}\frac{i_{F\wedge
_{E}B}^{E}}{A}\circ
p_{A}^{F\wedge _{E}B})\circ \overline{\triangle }_{F\wedge _{E}B} \\
&=&(\frac{E}{F}\cot _{E}\frac{i_{F\wedge _{E}B}^{E}}{A})\circ
(p_{F}^{E}\circ \J{F}{B}\cot _{E}p_{A}^{F\wedge _{E}B})\circ \overline{%
\triangle }_{F\wedge _{E}B}.
\end{eqnarray*}
Now, since $\frac{E}{F}\cot _{E}\frac{i_{F\wedge _{E}B}^{E}}{A}$
is a monomorphism, we conclude.
\end{proof}

\begin{proposition}\label{pro: naturality of eta}
  Let $\M$ be an abelian monoidal category. Consider the following commutative diagrams in $\M$:
\begin{equation*}
\begin{tabular}{cc}
$\xymatrix@C=2cm{
 A_1\ar[d]_a\ar@{^{(}->}[r]^{i_{A_1}^{B_1}}& B_1 \ar[d]_{b} \ar@{^{(}->}[r]^{i_{B_1}^{E_1}}
                & E_1 \ar[d]^{e}  \\
  A_2\ar@{^{(}->}[r]_{i_{A_2}^{B_2}}& B_2 \ar@{^{(}->}[r]_{i_{B_2}^{E_2}}
                & E_2             }
$& $ \xymatrix@C=2cm{
  F_1 \ar[d]_{f} \ar@{^{(}->}[r]^{i_{F_1}^{E_1}}
                & E_1 \ar[d]^{e}  \\
  F_2 \ar@{^{(}->}[r]_{i_{F_2}^{E_2}}
                & E_2             }
$
\end{tabular}
\end{equation*}where the morphisms are coalgebra homomorphisms. Then the
following diagram
\begin{equation*}
\xymatrix@C=3cm{
  F_1\wedge_{E_1} B_1 \ar[d]_{f\wedge_e b} \ar[r]^{\eta^{E_1}(F_1,B_1,A_1)}
                & \frac{E_1}{F_1} \cot_{E_1}\frac{B_1}{A_1}\ar[d]^{\frac{e}{f} \cot_{e}\frac{b}{a} }\\
  F_2\wedge_{E_2} B_2 \ar[r]_{\eta^{E_2}(F_2,B_2,A_2)}
                & \frac{E_2}{F_2} \cot_{E_2}\frac{B_2}{A_2}            }
\end{equation*}
  is commutative.
\end{proposition}

\begin{proof}
We have:
\begin{eqnarray*}
&&(\frac{E_{2}}{F_{2}}\square
_{E_{2}}\frac{i_{B_{2}}^{E_{2}}}{A_{2}})\circ
\eta ^{E_{2}}(F_{2},B_{2},A_{2})\circ (f\wedge _{e}b) \\
&=&(p_{F_{2}}^{E_{2}}\square _{E_{2}}p_{A_{2}}^{E_{2}})\circ \overline{%
\Delta }_{E_{2}}^{E_{2}}\circ i_{F_{2}\wedge
_{E_{2}}B_{2}}^{E_{2}}\circ
(f\wedge _{e}b) \\
&=&(p_{F_{2}}^{E_{2}}\square _{E_{2}}p_{A_{2}}^{E_{2}})\circ \overline{%
\Delta }_{E_{2}}^{E_{2}}\circ e\circ i_{F_{1}\wedge _{E_{1}}B_{1}}^{E_{1}} \\
&=&(p_{F_{2}}^{E_{2}}\circ e\square _{E_{2}}p_{A_{2}}^{E_{2}}\circ
e)\circ \overline{\Delta }_{E_{1}}^{E_{2}}\circ i_{F_{1}\wedge
_{E_{1}}B_{1}}^{E_{1}}
\\
&=&(\frac{e}{f}\circ p_{F_{1}}^{E_{1}}\square
_{E_{2}}\frac{e}{a}\circ p_{A_{1}}^{E_{1}})\circ \overline{\Delta
}_{E_{1}}^{E_{2}}\circ
i_{F_{1}\wedge _{E_{1}}B_{1}}^{E_{1}} \\
&=&(\frac{e}{f}\square _{E_{2}}\frac{e}{a})\circ
(p_{F_{1}}^{E_{1}}\square _{E_{2}}p_{A_{1}}^{E_{1}})\circ
\overline{\Delta }_{E_{1}}^{E_{2}}\circ
i_{F_{1}\wedge _{E_{1}}B_{1}}^{E_{1}} \\
&=&(\frac{e}{f}\square _{E_{2}}\frac{e}{a})\circ
(p_{F_{1}}^{E_{1}}\square
_{E_{2}}p_{A_{1}}^{E_{1}})\circ (E_{1}\square _{e}E_{1})\circ \overline{%
\Delta }_{E_{1}}^{E_{1}}\circ i_{F_{1}\wedge _{E_{1}}B_{1}}^{E_{1}} \\
&=&(\frac{e}{f}\square _{E_{2}}\frac{e}{a})\circ
(\frac{E_{1}}{F_{1}}\square _{e}\frac{E_{1}}{A_{1}})\circ
(p_{F_{1}}^{E_{1}}\square _{E_{1}}p_{A_{1}}^{E_{1}})\circ
\overline{\Delta }_{E_{1}}^{E_{1}}\circ
i_{F_{1}\wedge _{E_{1}}B_{1}}^{E_{1}} \\
&=&(\frac{e}{f}\square _{E_{2}}\frac{e}{a})\circ
(\frac{E_{1}}{F_{1}}\square
_{e}\frac{E_{1}}{A_{1}})\circ (\frac{E_{1}}{F_{1}}\square _{E_{1}}\frac{%
i_{B_{1}}^{E_{1}}}{A_{1}})\circ \eta ^{E_{1}}(F_{1},B_{1},A_{1}) \\
&=&(\frac{e}{f}\square _{e}\frac{e}{a}\circ \frac{i_{B_{1}}^{E_{1}}}{A_{1}}%
)\circ \eta ^{E_{1}}(F_{1},B_{1},A_{1}) \\
&=&(\frac{e}{f}\square _{e}\frac{i_{B_{2}}^{E_{2}}}{A_{2}}\circ \frac{b}{a}%
)\circ \eta ^{E_{1}}(F_{1},B_{1},A_{1}) \\
&=&(\frac{E_{2}}{F_{2}}\square
_{E_{2}}\frac{i_{B_{2}}^{E_{2}}}{A_{2}})\circ (\frac{e}{f}\square
_{e}\frac{b}{a})\circ \eta ^{E_{1}}(F_{1},B_{1},A_{1}).
\end{eqnarray*}%
Since $\frac{E_{2}}{F_{2}}\square
_{E_{2}}\frac{i_{B_{2}}^{E_{2}}}{A_{2}}$ is a monomorphism, we
finally obtain:
\begin{equation*}
\eta ^{E_{2}}(F_{2},B_{2},A_{2})\circ (f\wedge
_{e}b)=(\frac{e}{f}\square _{e}\frac{b}{a})\circ \eta
^{E_{1}}(F_{1},B_{1},A_{1}).
\end{equation*}
\end{proof}

\begin{proposition}\label{pro: diagram eta}
Let $\M$ be an abelian monoidal category. Let
$i_F=i_F^E:F\hookrightarrow E$, $i_A^B:A\hookrightarrow B$ and
$i_B=i_B^E:B\hookrightarrow E$ be monomorphisms which are
coalgebra homomorphisms. Then the following diagram
\begin{equation}\label{diagram eta}
\xymatrix@C=3cm{
  F\wedge B \ar[dr]_{\eta(F,B,A)} \ar[r]^{\eta(F,B,0)}
                & \frac{E}{F}\cot_E B \ar[d]^{\frac{E}{F}\cot_E p^B_A}  \\
                & \frac{E}{F}\cot_E \frac{B}{A}             }
\end{equation}is commutative. Furthermore, if the morphism $E/F \cot_E
p_A^B:{^E\M}\to \M$ is an epimorphism we get the following exact
sequence:
\begin{equation}\label{sequence eta epi}
\xymatrix@C=0.9cm{
  0 \ar[r] & F\w A \ar[rr]^{i_{F\w A}^{F\w B}} && F\w B\ar[rr]^{\eta(F,B,A)}&&\frac{E}{F}\cot_E \frac{B}{A} \ar[r]&0.}
\end{equation}
\end{proposition}

\begin{proof}
We apply Proposition \ref{pro: naturality of eta} in the case when
\begin{gather*}
E_{1}=E_{2}=E,\text{ }F_{1}=F_{2}=F,\text{ }B_{1}=B_{2}=B,\text{ }A_{1}=0,%
\text{ }A_{2}=A \\
e=\Id_E,\text{ }a=0,\text{ }f=\mathrm{Id}_{F},\text{
}b=\mathrm{Id}_{B}.
\end{gather*}
If the morphism $E/F\cot _{E}p_{A}^{B}$ is an epimorphism, then
$\eta (F,B,A)$ is an epimorphism as a composition of epimorphisms.
Thus, in view of (\ref{sequence eta}) we obtain (\ref{sequence eta
epi}).
\end{proof}

\begin{lemma}
\label{lem: F wedge B/F} Let $i_{F}:F\rightarrow E$ and
$i_{B}:B\rightarrow E $ be monomorphisms which are coalgebra
homomorphisms in an abelian monoidal category $\mathcal{M}$. Let
\begin{equation*}
(L,p):=\text{Coker}(i_{B}^{F\wedge _{E}B})=\frac{F\wedge
_{E}B}{B}.
\end{equation*}%
Then there is a unique morphism $^{F}\rho _{L}:L\rightarrow
F\otimes L$ such that
\begin{equation*}
^{E}\rho _{L}=(i_{F}\otimes L)\circ {^{F}\rho _{L}}.
\end{equation*}%
Moreover $(L,{^{F}\rho _{L}})$ is a left $F$-comodule and ${^{F}}\overline{{%
\rho }}{_{L}:L\rightarrow F}\square _{F}L$ is a morphism of left $E$%
-comodules.
\end{lemma}

\begin{proof}
Tensorize the following exact sequence
\begin{equation*}
0\rightarrow F\overset{i_{F}}{\longrightarrow }E\overset{p_{F}}{%
\longrightarrow }\frac{E}{F}\rightarrow 0
\end{equation*}%
by $L$ to get the exact sequence
\begin{equation}\label{formula ma1}
\xymatrix@C=0.5cm{
  0 \ar[r] & F\ot L \ar[rr]^{i_F\ot L} && E\ot L \ar[rr]^{p_F \ot L} && \frac{E}{F}\ot L.}
\end{equation}%
We have
\begin{eqnarray*}
&&(E/F\otimes i_{L}^{E/B})\circ (p_{F}\otimes L)\circ {^{E}\rho
_{L}}\circ p
\\
&=&(E/F\otimes i_{L}^{E/B})\circ (p_{F}\otimes L)\circ (E\otimes p)\circ {%
^{E}\rho _{F\wedge _{E}B}} \\
&=&(p_{F}\otimes i_{L}^{E/B}\circ p)\circ {^{E}\rho _{F\wedge _{E}B}} \\
&=&(p_{F}\otimes i_{L}^{E/B}\circ p)\circ \lbrack \J{F}{B}\otimes
(F\wedge
_{E}B)]\circ {\Delta }_{F\wedge _{E}B} \\
&=&\left( p_{F}\circ \J{F}{B}\otimes p_{B}\circ \J{F}{B}\right) \circ {\Delta }%
_{F\wedge _{E}B} \\
&=&\left( p_{F}\otimes p_{B}\right) \circ {\Delta }_{E}\circ
\J{F}{B}=0.
\end{eqnarray*}%
As $E/F\otimes i_{L}^{E/B}$ is a monomorphism and $p$ is an
epimorphism, we
obtain $(p_{F}\otimes L)\circ {^{E}\rho _{L}=0.}$ Since, by \ref{formula ma1}%
, we have
\begin{equation*}
(F\otimes L,i_{F}\otimes L)=\mathrm{Ker\,}(p_{F}\otimes L),
\end{equation*}%
there exists a unique morphisms $^{F}\rho _{L}:L\rightarrow
F\otimes L$ such that $^{E}\rho _{L}=(i_{F}\otimes L)\circ
{^{F}\rho _{L}}.$ Moreover, by Lemma \ref{lem: induced comodule}
$(L,^{F}\rho _{L})$ is a left $F$-comodule and
${^{F}}\overline{{\rho }}{_{L}}:L\rightarrow F\square _{F}L$ is a
morphism of left $E$-comodules.
\end{proof}

\begin{lemma}
\label{lem: rho inscatolati} Let $\alpha :F\rightarrow A$ be a
homomorphism and let $i_{A}^{E}:A\rightarrow E$ be a monomorphism
which is a coalgebra homomorphism in an abelian monoidal category
$\mathcal{M}$. Let $(W,^{A}\rho _{W})$ be a left $A$-comodule. For
a morphism $^{F}\rho _{W}:W\rightarrow F\otimes W$ the following
statement are equivalent.

(1) $^{A}\rho _{W}=(\alpha \otimes W)\circ {^{F}\rho _{W}.}$

(2) $^{E}\rho _{W}=(i_{A}^{E}\alpha \otimes W)\circ {^{F}\rho
}{_{W}}$, where $^{E}\rho _{W}=(i_{A}^{E}\otimes W)\circ {^{A}\rho
_{W}.}$
\end{lemma}

\begin{proof}
$\left( 1\right) \Rightarrow \left( 2\right) $ is trivial.\newline
$\left( 2\right) \Rightarrow \left( 1\right) $ We have that $%
(i_{A}^{E}\otimes W)\circ {^{A}\rho _{W}}={^{E}\rho
_{W}}=(i_{A}^{E}\alpha \otimes W)\circ {^{F}\rho
}{_{W}}=(i_{A}^{E}\otimes W)\circ (\alpha \otimes W)\circ
{^{F}\rho }{_{W}}$. Since $i_{A}^{E}\otimes W$ is a monomorphism,
we have $^{A}\rho _{W}=(\alpha \otimes W)\circ {^{F}\rho }{_{W}.}$
\end{proof}

\begin{proposition}
\label{pro: rho chiave} Let $\alpha :F\rightarrow A$ be a
coalgebra homomorphism, let $i_{B}^{A}:B\rightarrow A$ and
$i_{A}^{E}:A\rightarrow E$ be monomorphisms which are coalgebra
homomorphisms in an abelian monoidal category $\mathcal{M}$.
Assume there is a homomorphism $\pi :A\rightarrow F$ such that
$\pi \alpha =\mathrm{Id}_{F}.$ Assume that there exists a morphism
$^{F}\rho _{A/B}:A/B\rightarrow F\otimes A/B$ such that $^{E}\rho
_{A/B}=(i_{A}^{E}\alpha \otimes A/B)\circ {^{F}\rho _{A/B},}$ where $%
^{E}\rho _{A/B}=(i_{A}^{E}\otimes A/B)\circ {^{A}\rho
_{A/B}.}$\newline Then we have
\begin{equation*}
^{F}\rho _{A/B}=(\pi \otimes A/B)\circ {^{A}\rho _{A/B}.}
\end{equation*}%
Furthermore ${^{F}\rho _{A/B}}$ is uniquely defined by the
following relation
\begin{equation*}
{^{F}\rho _{A/B}}\circ p_{B}^{A}=(\pi \otimes p_{B}^{A})\circ
{\Delta }_{A}
\end{equation*}%
and we have
\begin{equation*}
(A\otimes p_{B}^{A})\circ {\Delta }_{A}=(\alpha \pi \otimes p_{B}^{A})\circ {%
\Delta }_{A}.
\end{equation*}
\end{proposition}

\begin{proof}
By Lemma \ref{lem: rho inscatolati} one has $^{A}\rho
_{A/B}=(\alpha \otimes A/B)\circ {^{F}\rho _{A/B},}$ so that we
have
\begin{equation*}
^{F}\rho _{A/B}=(\pi \otimes A/B)\circ (\alpha \otimes A/B)\circ
{^{F}\rho _{A/B}=}(\pi \otimes A/B)\circ {^{A}\rho _{A/B}}.
\end{equation*}%
Furthermore we get:
\begin{equation*}
{^{F}\rho _{A/B}}\circ p_{B}^{A}=(\pi \otimes A/B)\circ {^{A}\rho
_{A/B}}\circ p_{B}^{A}=(\pi \otimes A/B)\circ (A\otimes p_{B}^{A})\circ {%
\Delta }_{A}=(\pi \otimes p_{B}^{A})\circ {\Delta }_{A}{.}
\end{equation*}%
Finally we obtain:
\begin{eqnarray*}
(\alpha \pi \otimes p_{B}^{A})\circ {\Delta }_{A} &=&(\alpha
\otimes
A/B)\circ (\pi \otimes p_{B}^{A})\circ {\Delta }_{A} \\
&=&(\alpha \otimes A/B)\circ {^{F}\rho _{A/B}\circ }p_{B}^{A} \\
&=&{^{A}\rho _{A/B}}\circ p_{B}^{A}=(A\otimes p_{B}^{A})\circ
{\Delta }_{A}.
\end{eqnarray*}
\end{proof}

\begin{corollary}
Let $i_{F}:F\rightarrow E$ and $i_{B}:B\rightarrow E$ be
monomorphisms which are coalgebra homomorphisms in an abelian
monoidal category $\mathcal{M}$. Let
\begin{equation*}
(L,p):=\text{Coker}(i_{B}^{F\wedge _{E}B})=\frac{F\wedge
_{E}B}{B}.
\end{equation*}%
Assume that $\pi :F\wedge _{E}B\rightarrow F$ is a morphism such
that $\pi \circ i_{F}^{F\wedge _{E}B}=\mathrm{Id}_{F}.$ Then the
morphism $^{F}\rho _{L}:L\rightarrow F\otimes L$ defined in Lemma
\ref{lem: F wedge B/F} is uniquely defined by the following
relation
\begin{equation}
{^{F}\rho _{L}}\circ p_{B}^{F\wedge _{E}B}=(\pi \otimes
p_{B}^{F\wedge _{E}B})\circ {\Delta }_{F\wedge _{E}B}
\label{formula: eta is comodule morph}
\end{equation}%
(which means that $p_{B}^{F\wedge _{E}B}$ is a morphism of left $F$%
-comodules whenever $\pi$ is a coalgebra homomorphism) and we
have:
\begin{equation}
[(F\wedge _{E}B)\otimes p_{B}^{F\wedge _{E}B}]\circ {\Delta
}_{F\wedge _{E}B}=(i_{F}^{F\wedge _{E}B}\pi \otimes p_{B}^{F\wedge
_{E}B})\circ {\Delta }_{F\wedge _{E}B}.  \label{formula: f_D
disappears}
\end{equation}
\end{corollary}

\begin{proof}
Since, by Lemma \ref{lem: F wedge B/F}, one has
\begin{equation*}
^{E}\rho _{L}=(i_{F}\otimes L)\circ {^{F}\rho _{L}}=(\J{F}{B}\circ
i_{F}^{F\wedge _{E}B}\otimes L)\circ {^{F}\rho _{L}}
\end{equation*}%
so that, we can apply Proposition \ref{pro: rho chiave} in the case when $%
A=F\wedge _{E}B,B=B,E=E,F=F,\alpha =i_{F}^{F\wedge
_{E}B},i_{A}^{E}=\J{F}{B}.$
\end{proof}

\begin{claim}
\label{claim: golden claim}Let $\alpha :F\rightarrow A$ be a
coalgebra homomorphism and let $i_{A}^{E}:A\rightarrow E$ be a
monomorphism which is a
coalgebra homomorphism in an abelian monoidal category $\mathcal{M}$. Let $%
\;(L,p)=\text{Coker}(\alpha ).$ Since $F$ and $A$ are left
$E$-comodules via
$i_{A}^{E}\alpha $ and $i_{A}^{E}$ respectively, so is $L.$ Its left $E$%
-comodule structure is uniquely defined by a morphism $^{E}\rho
_{L}:L\rightarrow E\otimes L$ such that
\begin{equation*}
^{E}\rho _{L}\circ p=(E\otimes p)\circ {^{E}\rho
_{A}}=(i_{A}^{E}\otimes p)\circ {\Delta }_{A}.
\end{equation*}%
We also point out that, since $F$ and $A$ are left $A$-comodules
via $\alpha $ and $\mathrm{Id}_{A}$ respectively, so is $L.$ Its
left $A$-comodule structure is uniquely defined by a morphism
$^{A}\rho _{L}:L\rightarrow A\otimes L$ such that
\begin{equation*}
^{A}\rho _{L}\circ p=(A\otimes p)\circ {^{A}\rho _{A}}=(A\otimes p)\circ {%
\Delta }_{A}.
\end{equation*}%
\newline
Assume that there exists a morphism
\begin{equation*}
^{F}\rho _{L}:L\rightarrow F\otimes L
\end{equation*}%
such that $^{E}\rho _{L}=(i_{A}^{E}\alpha \otimes L)\circ
{^{F}\rho _{L}}$ (this happens, for example, in the case
$A=F\wedge _{E}F$ and $\alpha
=i_{F}^{F\wedge _{E}F}$). Then, by Lemma \ref{lem: induced comodule}, $%
(L,^{F}\rho _{L})$ is a left $F$-comodule. Hence one can endow $L$
with a left $A$-comodule structure via $^{A}\rho _{L}^{\prime
}=(\alpha \otimes L)\circ ^{F}\rho _{L}.$ One has
\begin{equation*}
^{A}\rho _{L}^{\prime }={^{A}\rho _{L}}.
\end{equation*}%
In fact we have
\begin{eqnarray*}
(i_{A}^{E}\otimes L)\circ {^{A}\rho _{L}^{\prime }\circ p}
&=&(i_{A}^{E}\alpha
\otimes L)\circ {^{F}\rho _{L}\circ p} \\
&=&^{E}\rho _{L}\circ p=(i_{A}^{E}\otimes p)\circ {\Delta }%
_{A}\\&=&(i_{A}^{E}\otimes L)\circ (A\otimes p)\circ {\Delta }%
_{A}=(i_{A}^{E}\otimes L)\circ {^{A}\rho _{L}\circ p}.
\end{eqnarray*}%
Since $i_{A}^{E}\otimes L$ is a monomorphism and $p$ is an
epimorphism, we conclude.\newline
Assume there is a homomorphism $\pi :A\rightarrow F$ such that $\pi \alpha =%
\mathrm{Id}_{F}.$ Then we have
\begin{equation*}
(F\otimes p)\circ (\pi \otimes A)\circ {\Delta }_{A}\circ \alpha
=(\pi \otimes p)\circ (\alpha \otimes \alpha )\circ {\Delta
}_{F}=0.
\end{equation*}%
Since $(L,p)=\text{Coker}(\alpha ),$ there exists a morphism
$^{F}\rho _{L}^{\prime \prime }:L\rightarrow F\otimes L$ such that
\begin{equation*}
{^{F}\rho _{L}^{\prime \prime }\circ p}=(\pi \otimes {p})\circ
{\Delta }_{A}
\end{equation*}%
(when $\pi $ is a coalgebra morphism, ${^{F}\rho _{L}^{\prime
\prime }}$
defines the left $F$-comodule structure that $L$ has via $\pi )$. By Lemma %
\ref{pro: rho chiave}, we have that ${^{F}\rho _{L}}$ is uniquely
defined by the following relation
\begin{equation*}
{^{F}\rho _{L}\circ p}=(\pi \otimes {p})\circ {\Delta }_{A}.
\end{equation*}%
Therefore
\begin{equation*}
^{F}\rho _{L}^{\prime \prime }={^{F}\rho _{L}}.
\end{equation*}
\end{claim}

\begin{proposition}
\label{lem: passo al bar}Let $\alpha :F\rightarrow A$ be a
coalgebra homomorphism and let $i_{A}^{E}:A\rightarrow E$ be a
monomorphism which is a coalgebra homomorphism in an abelian
monoidal category $\mathcal{M}$. For any left $F$-comodule
$(W,^{F}\rho _{W})$, let
\begin{equation*}
^{A}\rho _{W}=(\alpha \otimes W)^{F}\rho _{W},\text{\qquad
}^{E}\rho _{W}=(i_{A}^{E}\alpha \otimes W)^{F}\rho _{W}\text{.}
\end{equation*}%
Let $f:W_{1}\rightarrow W_{2}$ be a morphism of left
$F$-comodules. Then
\begin{equation*}
(i_{A}^{E}\square _{F}f)\circ ^{A}\overline{\rho }_{W_{1}}=^{E}\overline{%
\rho }_{W_{2}}\circ f.
\end{equation*}
\end{proposition}

\begin{proof}
We have
\begin{equation*}
^{E}\rho _{W_{2}}\circ f=(i_{A}^{E}\alpha \otimes W_{2})\circ
{^{F}\rho _{W_{2}}}\circ f=(i_{A}^{E}\alpha \otimes W_{2})\circ
(F\otimes f)\circ {^{F}\rho _{W_{1}}}=(i_{A}^{E}\otimes f)\circ
{^{A}\rho _{W_{1}}},
\end{equation*}%
so that
\begin{eqnarray*}
&&\chi _{F}(E,W_{2})\circ {^{E}\overline{\rho}}_{W_{2}}\circ f \\
&=&^{E}\rho _{W_{2}}\circ f\\
&=&(i_{A}^{E}\otimes f)\circ {^{A}\rho
_{W_{1}}}=(i_{A}^{E}\otimes f)\circ \chi _{F}(A,W_{1})\circ {^{A}\overline{%
\rho }_{W_{1}}}=\chi _{F}(E,W_{2})\circ (i_{A}^{E}\square _{F}f)\circ {^{A}%
\overline{\rho }_{W_{1}}}.
\end{eqnarray*}%
Since $\chi _{F}(E,W_{2})$ is a monomorphism, we conclude.
\end{proof}

\section{Main results}
We now recall some definitions and results established in
\cite{Cotensor}.
\begin{claim}
\label{def of delta_n}Let $X$ be an object in an abelian monoidal category $%
\left( \mathcal{M},\otimes ,\mathbf{1}\right) $. Set%
\begin{equation*}
X^{\otimes 0}=\mathbf{1},\qquad X^{\otimes 1}=X\qquad
\text{and}\qquad X^{\otimes n}=X^{\otimes n-1}\otimes X,\text{ for
every }n>1
\end{equation*}%
and for every morphism $f:X\rightarrow Y$ in $\mathcal{M}$, set%
\begin{equation*}
f^{\otimes 0}=\mathrm{Id}_{\mathbf{1}},\qquad f^{\otimes 1}=f\qquad \text{and%
}\qquad f^{\otimes n}=f^{\otimes n-1}\otimes f,\text{ for every
}n>1.
\end{equation*}%
Let $\left( C,\Delta _{C},\varepsilon _{C}\right) $ be a coalgebra in $%
\mathcal{M}$ and for every $n\in
\mathbb{N}
,$ define the $n^{\text{th}}$ iterated comultiplication of $C,$
$\Delta _{C}^{n}:C\rightarrow C^{\otimes {n+1}}$, by
\begin{equation*}
\Delta _{C}^{0}=\text{Id}_{C},\qquad \Delta _{C}^{1}=\Delta _{C}\qquad \text{%
and}\qquad \Delta _{C}^{n}=\left( \Delta _{C}^{\otimes n-1}\otimes
C\right) \Delta _{C},\text{ for every }n>1.
\end{equation*}%
Let $\delta :D\rightarrow E$ be a monomorphism which is a
homomorphism of
coalgebras in $\mathcal{M}$. Denote by $(L,p)$ the cokernel of $\delta $ in $%
\mathcal{M}$. Regard $D$ as a $E$-bicomodule via $\delta $ and
observe that $L$ is a $E$-bicomodule and $p$ is a morphism of
bicomodules. Let
\begin{equation*}
(D^{\wedge _{E}^{n}},\delta _{n}):=\ker (p^{\otimes {n}}\Delta
_{E}^{n-1})
\end{equation*}%
for any $n\in \mathbb{N}\setminus \{0\}.$ Note that $(D^{\wedge
_{E}^{1}},\delta _{1})=(D,\delta )$ and $(D^{\wedge
_{E}^{2}},\delta _{2})=D\wedge _{E}D.$ \newline In order to
simplify the notations we set $(D^{\wedge _{E}^{0}},\delta
_{0})=(0,0).$\newline Now, since $\mathcal{M}$ has left exact
tensor functors and since $p^{\otimes {n}}\Delta _{E}^{n-1}$ is a
morphism of $E$-bicomodules (as a composition of morphisms
of $E$-bicomodules), we get that $D^{\wedge _{E}^{n}}$ is a coalgebra and $%
\delta _{n}:D^{\wedge _{E}^{n}}\rightarrow E$ is a coalgebra
homomorphism for any $n>0$ and hence for any $n\in \mathbb{N}$.
\end{claim}

\begin{proposition}\cite[Proposition 1.10]{Cotensor}\label{pro: limit of delta}Let $\delta:D\to E$ be a monomorphism which is a morphism of coalgebras in an abelian monoidal category
$\M$. Then for any $i\leq j$ in $\mathbb{N}$ there is a (unique)
morphism $\xi_{i}^j:D^{\wedge_E ^i}\to D^{\wedge_E ^j}$ such that
\begin{equation}\label{compatibility of delta i}
\delta_j\xi_{i}^j=\delta_i.
\end{equation} Moreover $\xi_{i}^j$ is a coalgebra
homomorphism and $((D^{\wedge_E ^i})_{i\in \mathbb{N}},(\xi
_{i}^j)_{i,j\in \mathbb{N}})$ is a direct system in $\M$ whose
direct limit, if it exists, carries a natural coalgebra structure
that makes it the direct limit of $((D^{\wedge_E ^i})_{i\in
\mathbb{N}},(\xi _{i}^j)_{i,j\in \mathbb{N}})$ as a direct system
of coalgebras.
\end{proposition}

\begin{notation}\label{notation tilde}Let $\delta:D\to E$ be a morphism of coalgebras in an abelian monoidal category
$\M$ cocomplete and with left exact tensor functors . By
Proposition \ref{pro: limit of delta} $((D^{\wedge_E ^i})_{i\in
\mathbb{N}},(\xi _{i}^j)_{i,j\in \mathbb{N}})$ is a direct system
in $\M$ whose direct limit carries a natural coalgebra structure
that makes it the direct limit of $((D^{\wedge_E ^i})_{i\in
\mathbb{N}},(\xi _{i}^j)_{i,j\in \mathbb{N}})$ as a direct system
of coalgebras.\\From now on we set: $(\widetilde{D}_E,
(\xi_i)_{i\in \mathbb{N}})= \underrightarrow{\lim }(D^{\wedge_E
^i})_{i\in \mathbb{N}}$, where $\xi_i:D^{\wedge_E ^i}\to
\widetilde{D}_E$ denotes the structural morphism of the direct
limit. We simply write $\widetilde{D}$ if there is no danger of
confusion. We note that, since $\widetilde{D}$ is a direct limit
of coalgebras, the canonical (coalgebra) homomorphisms $(\delta
_{i}:D^{\wedge_E ^i}\rightarrow E)_{i\in \mathbb{N}}$, which are
compatible by (\ref{compatibility of delta i}), factorize to a
unique coalgebra homomorphism $\widetilde{\delta}
:\widetilde{D}\rightarrow E$ such that $\widetilde{\delta}
\xi_i=\delta_i$ for any $i\in \mathbb{N}$.
\end{notation}

\begin{claim}
\label{claim 4.2}Let $(C,\Delta ,\varepsilon )$ be a coalgebra in
a cocomplete abelian monoidal category $\mathcal{M}$ and let $(M,\rho _{M}^{r},\rho _{M}^{l})$ be a $C$--$\,$%
{}{}bicomodule. Set
\begin{equation*}
M^{\square 0}=C,M^{\square 1}=M\text{\quad and\quad }M^{\square
n}=M^{\square n-1}\square M\text{ for any }n>1
\end{equation*}%
and define $(C^{n}(M))_{n\in \mathbb{N}}$ by
\begin{equation*}
C^{0}(M)=0,C^{1}(M)=C\text{\quad and\quad
}C^{n}(M)=C^{n-1}(M)\oplus M^{\square n-1}\text{ for any }n>1.
\end{equation*}%
Let $\sigma _{i}^{i+1}:C^{i}(M)\rightarrow C^{i+1}(M)$ be the
canonical inclusion and for any $j>i,$ define:
\begin{equation*}
\sigma _{i}^{j}=\sigma _{j-1}^{j}\sigma _{j-2}^{j-1}\cdots \sigma
_{i+1}^{i+2}\sigma _{i}^{i+1}:C^{i}(M)\rightarrow C^{j}(M).
\end{equation*}%
Then $(({C^{i}}(M))_{i\in \mathbb{N}},(\sigma _{i}^{j})_{i,j\in
\mathbb{N}})$ is a direct system in $\mathcal{M}$. We set
\begin{equation*}
T_{C}^{c}(M)=\bigoplus_{n\in \mathbb{N}}M^{\square n}=C\oplus
M\oplus M^{\square 2}\oplus M^{\square 3}\oplus \cdots
\end{equation*}%
and we denote by $\sigma _{i}:C^{i}(M)\rightarrow T_{C}^{c}(M)$
the canonical inclusion.\medskip\newline Throughout let
\begin{eqnarray*}
\pi _{n}^{m} &:&C^{n}(M)\rightarrow C^{m}(M)\text{ }(m\leq n),\text{\qquad }%
\pi _{n}:T_{C}^{c}(M)\rightarrow C^{n}(M), \\
p_{n}^{m} &:&C^{n}(M)\rightarrow M^{\square m}\text{ }(m<n),\text{\qquad }%
p_{n}:T_{C}^{c}(M)\rightarrow M^{\square n},
\end{eqnarray*}
be the canonical projections and let
\begin{eqnarray*}
\sigma _{m}^{n} &:&C^{m}(M)\rightarrow C^{n}(M)\text{ }(m\leq n),\text{%
\qquad }\sigma _{n}:C^{n}(M)\rightarrow T_{C}^{c}(M), \\
i_{m}^{n} &:&M^{\square m}\rightarrow C^{n}(M)\text{ }(m<n),\text{\qquad }%
i_{m}:M^{\square m}\rightarrow T_{C}^{c}(M),
\end{eqnarray*}%
be the canonical injection for any $m,n\in \mathbb{N}$.\newline
For technical reasons we set $\pi _{n}^{m}=0,$ $\sigma _{m}^{n}=0$ for any $n<m$ and $%
p_{n}^{m}=0$, $i_{m}^{n}=0$ for any $n\leq m.$ Then, we have the
following relations:%
\begin{equation*}
p_{n}\sigma _{k}=p_{k}^{n},\qquad p_{n}i_{k}=\delta _{n,k}\mathrm{Id}%
_{M^{\square k}},\qquad \pi _{n}i_{k}=i_{k}^{n}.
\end{equation*}%
Moreover, we have:
\begin{equation*}
\begin{tabular}{lll}
$\pi _{n}^{m}\sigma _{k}^{n}=\sigma _{k}^{m},\text{ if }k\leq m\leq n,$ & $%
\qquad $and$\qquad $ & $\pi _{n}^{m}\sigma _{k}^{n}=\pi _{k}^{m},\text{ if }%
m\leq k\leq n,$ \\
$p_{n}^{m}\pi _{k}^{n}=p_{k}^{m},\text{ if }m<n\leq k,$ & $\qquad $and$%
\qquad $ & $\sigma _{n}^{m}i_{k}^{n}=i_{k}^{m},\text{ if }k<n\leq m,$ \\
$p_{n}^{m}\sigma _{k}^{n}=p_{k}^{m},\text{ if }m<k\leq n,$ & $\qquad $and$%
\qquad $ & $\pi _{n}^{m}i_{k}^{n}=i_{k}^{m},\text{ if }k<m\leq n,$ \\
$p_{n}^{m}\pi _{n}=p_{m},\text{ if }m<n,$ & $\qquad $and$\qquad $
& $\sigma
_{n}i_{m}^{n}=i_{m},\text{ if }m<n\text{,}$ \\
$\pi _{n}\sigma _{k}=\sigma _{k}^{n},\text{ if }k\leq n,$ & $\qquad $and$%
\qquad $ & $\pi _{n}\sigma _{k}=\pi _{k}^{n},\text{ if }n\leq k,$ \\
$p_{n}^{m}i_{m}^{n}=\mathrm{Id}_{M^{\square m}},\text{ if }m<n.$ &
&
\end{tabular}%
\end{equation*}
In the other cases, these compositions are zero.
\end{claim}

Let $(C,\Delta ,\varepsilon )$ be a coalgebra in a cocomplete
abelian monoidal category $\M$ with
left exact tensor functors and let $(M,\rho _{M}^{r},\rho _{M}^{l})$ be a $%
C$-bicomodule.
Then $ (T^c_C(M),(\sigma_n)_{n\in \N})=\underrightarrow{%
\lim }C^{i}(M).$

\begin{theorem}\cite[Theorem 2.9]{Cotensor}\label{teo: T as limit}Let $(C,\Delta ,\varepsilon )$ be a
coalgebra in a cocomplete abelian monoidal category $\M$ and let
$(M,\rho _{M}^{r},\rho _{M}^{l})$ be a $C$-bicomodule.
  $ (T^c_C(M),(\sigma_i)_{i\in \N})$ carries a
natural coalgebra structure that makes it the direct limit of
$((C^{i}(M))_{i\in \mathbb{N}},(\sigma _{i}^j)_{i,j\in
\mathbb{N}})$ as a direct system of coalgebras.
\end{theorem}

\begin{theorem}\cite[Theorem 2.13]{Cotensor}
\label{teo: pre univ property of cotensor coalgebra} Let
$(C,\Delta ,\varepsilon )$ be a coalgebra in a cocomplete abelian
monoidal category $\M$ and let $(M,\rho _{M}^{r},\rho _{M}^{l})$
be a $C$-bicomodule. Let $\delta :D\rightarrow E$ be a
monomorphism
which is a morphism of coalgebras such that the canonical morphism $%
\widetilde{\delta }:\widetilde{D}\rightarrow E$ of Notation
\ref{notation tilde} is a monomorphism. Let
$f_{C}:\widetilde{D}\rightarrow C$ be a coalgebra homomorphism and
let $f_{M}:\widetilde{D}\rightarrow M$ be a
morphism of $C$-bicomodules such that $f_{M}\xi _{1}=0$, where $%
\widetilde{D}$ is a $C$-bicomodule via $f_{C}.$ Then there is a
unique morphism $f:\widetilde{D}\rightarrow T_{C}^{c}(M)$ such
that
\begin{equation*}
f\xi _{n}=\sigma _{n}f_{n},\text{ for any }n\in \N,
\end{equation*}
where
\begin{equation}
f_{n}=\sum_{t=0}^{n}i_{t}^{n}f_{M}^{\square t}\overline{\Delta }_{\widetilde{%
D}}^{t-1}\xi _{n}  \label{formula: f_n}
\end{equation}
and $\overline{\Delta
}_{\widetilde{D}}^{n}:\widetilde{D}\rightarrow
\widetilde{D}^{\square n+1}$ is the $n^{\text{th}}$ iteration of $\overline{%
\Delta }_{\widetilde{D}}$ $(\overline{\Delta }_{\widetilde{D}}^{-1}=f_{C},%
\overline{\Delta }_{\widetilde{D}}^{0}=\Id_{\widetilde{D}},\overline{\Delta }%
_{\widetilde{D}}^{1}=\overline{\Delta }_{\widetilde{D}}:\widetilde{D}%
\rightarrow \widetilde{D}\square \widetilde{D}).$ \newline
Moreover:

1) $f$ is a coalgebra homomorphism;

2) $p_{0}f=f_{C}$ and $p_{1}f=f_{M}$, where
$p_{n}:T_{C}^{c}(M)\rightarrow M^{\square n}$ denotes the
canonical projection. \newline Furthermore, any coalgebra
homomorphism $f:\widetilde{D}\rightarrow T_{C}^{c}(M)$ that
fulfils 2) satisfies the following relation:
\begin{equation}
p_{k}f=f_{M}^{\square k}\overline{\Delta
}_{\widetilde{D}}^{k-1}\text{ for any }k\in \N.\label{relation f}
\end{equation}
\end{theorem}

\begin{theorem}\cite[Theorem 2.15]{Cotensor}
\label{coro: univ property of cotensor coalgebra} Let $(C,\Delta
,\varepsilon )$ be a coalgebra in a cocomplete and complete
abelian monoidal category $\M$ satisfying AB5. Let $(M,\rho
_{M}^{r},\rho _{M}^{l})$ be a $C$-bicomodule. Let $\delta
:D\rightarrow E$ be a monomorphism which is a homomorphism of coalgebras. Let $%
f_{C}:\widetilde{D}\rightarrow C$ be a coalgebra homomorphism and let $f_{M}:%
\widetilde{D}\rightarrow M$ be a morphism of $C$-bicomodules such
that $f_{M}\xi _{1}=0$, where $\widetilde{D}$ is a bicomodule via
$f_{C}.$ Then there is a unique coalgebra homomorphism
$f:\widetilde{D}\rightarrow
T_{C}^{c}(M)$ such that $p_{0}f=f_{C}$ and $p_{1}f=f_{M}$, where $%
p_{n}:T_{C}^{c}(M)\rightarrow M^{\square n}$ denotes the canonical
projection.
\begin{equation*}
\xymatrix@C=2cm{
  T_{C}^{c}(M) \ar[d]_{p_0} \ar[r]^{p_1} & M   \\
  C& \widetilde{D}  \ar[u]_{f_M}\ar@{.>}[ul]|-{f}\ar[l]^{f_C} &D\ar[l]^{\xi_1}  \ar[ul]_0          }
\end{equation*}
\end{theorem}

\begin{claim}\label{cl:E-Proj}
Let $\mathcal{M}$ be an abelian category and let $\mathcal{H}$ be
a class of monomorphisms in $\mathcal{M}$. We recall that an
object $I$ in $\mathcal{M}$ is called injective rel $\lambda $,
where $\lambda :X\rightarrow Y$ is a monomorphism in
$\mathcal{H}$, if $\mathcal{M}(\lambda,I ):\mathcal{
M}(Y,I)\rightarrow \mathcal{M}(X,I)$ is surjective. $I$ is called
$\mathcal{H}$-injective if it is injective rel $\lambda $ for
every $\lambda $ in $\mathcal{H}$. The\emph{\ closure} of
$\mathcal{H}$ is the class $ \mathcal{C(H)}$ containing all
monomorphisms $\lambda$ in $\mathcal{M}$ such that every
$\mathcal{H}$-injective object is also injective rel $ \lambda$.
The class $\mathcal{H}$ is called \emph{closed} if $\mathcal{H}$
is $\mathcal{C(H)}$. A closed class $\mathcal{H}$ is called
\emph{injective} if for any object $M$ in $\mathcal{M}$ there is
an monomorphism $\lambda:M\rightarrow I$ in $\mathcal{H}$ such
that $I$ is $\mathcal{H}$-injective.
\end{claim}

\begin{claim}\label{claim injective class}We fix a coalgebra $C$ in a monoidal category $\M$. Let $\mathbb{U}:\Mt\to\M$ be the forgetful functor. Then
\begin{equation} \mathcal{I}:=\{f\in
\Mt\mid \mathbb{U}(f)\text{ cosplits in }\M\}. \label{inj class}
\end{equation}
is an injective class of monomorphisms.\\ Now, for any
$C$-bicomodule $M\in\cmc$, we define the Hochschild cohomology of
$C$ with coefficients in $M$ by:
\[
\Hu(M,C)=\mathbf{Ext}^\bullet_{\mathcal{I}}(M,C) ,
\]
where $\mathbf{Ext}^\bullet_{\mathcal{I}}(M,-)$ are the relative
left derived functors of $\cmc(M,-)$. The notion of Hochschild
cohomology for algebras and coalgebras in monoidal categories has
been deeply investigated in \cite{AMS}. Here we quote some results
that will be needed afterwards.
\end{claim}

\begin{theorem}\cite[Theorem 4.22]{AMS}\label{X cof teo} \label{X coformally smooth teo}
\label{X formally smooth direct limit} Let $(D,\Delta ,\varepsilon
)$ be a coalgebra in an abelian monoidal category $\mathcal{M}$.
Then the following conditions are equivalent:

(a) $D$ is formally smooth (i.e. $\mathrm{H}^{2}\left( M,D\right) =0$, for any $M\in {}^{D}{{\mathcal{M}}}%
^{D}{.}$).

(b) The canonical morphism $\xi_1:D\rightarrow \widetilde{D}_E$
has a coalgebra homomorphism retraction, whenever:

\begin{itemize}
    \item $E$ is a coalgebra
endowed with a coalgebra homomorphism $\delta:D\rightarrow E$;
    \item $\delta$ is a monomorphism;
    \item $\widetilde{D}_E$ exists;
    \item for any $r\in \mathbb{N},$ the canonical injection $\xi
_{r}^{r+1}:D^{r}\rightarrow D^{r+1}$ cosplits in $\M.$
\end{itemize}
\end{theorem}

\begin{theorem}\cite[Theorem 4.15]{Cotensor}
\label{teo Cosmooth} Let $(C,\Delta ,\varepsilon )$ be a formally
smooth coalgebra in a cocomplete and complete abelian monoidal
category $\M$ satisfying $AB5,$ with left and right exact tensor
functors. Assume that denumerable coproducts commute with $\otimes
$. Let $(M,\rho _{M}^{r},\rho _{M}^{l})$ be a
$\mathcal{I}$-injective $C$-bicomodule. Then the cotensor
coalgebra $T_{C}^{c}(M)$ is formally smooth.
\end{theorem}

From now on we will use the following notation
\begin{equation*}
D^{n}:=D^{\w n},\text{ for every }n\in \mathbb{N.}
\end{equation*}

\begin{lemma}
\label{lem: lambda family}Let $\delta :D\rightarrow E$ be a
monomorphism which is a morphism of coalgebras in an abelian
monoidal category $\M$ and assume that for any $r\in \mathbb{N},$
the canonical injection $\xi _{r}^{r+1}:D^{r}\rightarrow D^{r+1}$
cosplits in $\M
$ i.e. there exists $\lambda _{r+1}^{r}:D^{r+1}\rightarrow D^{r}$ such that $%
\lambda _{r+1}^{r}\circ $ $\xi _{r}^{r+1}=\mathrm{Id}_{D^{r}}.$
Then, for any $r\in \mathbb{N},$ the canonical injection $\xi
_{r}:D^{r}\rightarrow \widetilde{D}$ cosplits in $\M.$
\end{lemma}

\begin{proof}
For every $i\in \mathbb{N}$ let us define $\lambda
_{i}:D^{i}\rightarrow D^{r}$ by setting
\begin{equation*}
\left\{
\begin{tabular}{l}
$\lambda _{i}=\lambda _{r+1}^{r}\circ \cdots \circ \lambda
_{i}^{i-1},\text{
if }i>r$ \\
$\lambda _{r}=\mathrm{Id}_{D^{r}}$ \\
$\lambda _{i}=\xi _{i}^{r}\text{ if }i<r.$%
\end{tabular}
\right.
\end{equation*}
Let us prove that $(\lambda _{i}:D^{i}\rightarrow D^{r})_{i\in
\mathbb{N}}$ is a compatible family of morphisms in $\M$ i.e. that
$\lambda _{i+1}\circ \xi _{i}^{i+1}=\lambda _{i}.$ \newline If
$i+1<r,$ we have
\begin{equation*}
\lambda _{i+1}\circ \xi _{i}^{i+1}=\xi _{i+1}^{r}\circ \xi
_{i}^{i+1}=\xi _{i}^{r}=\lambda _{i}.
\end{equation*}
If $i+1=r,$ we have
\begin{equation*}
\lambda _{i+1}\circ \xi _{i}^{i+1}=\mathrm{Id}_{D^{r}}\circ \xi
_{i}^{i+1}=\xi _{i}^{r}=\lambda _{i}.
\end{equation*}
If $i+1>r,$ we have
\begin{equation*}
\lambda _{i+1}\circ \xi _{i}^{i+1}=\lambda _{r+1}^{r}\circ \cdots
\circ \lambda _{i}^{i-1}\circ \lambda _{i+1}^{i}\circ \xi
_{i}^{i+1}=\lambda _{r+1}^{r}\circ \cdots \circ \lambda
_{i}^{i-1}=\lambda _{i}.
\end{equation*}
Now, since $(\lambda _{i}:D^{i}\rightarrow D^{r})_{i\in
\mathbb{N}}$ is a
compatible family of morphisms in $\M$ there exists a morphism $\lambda :%
\widetilde{D}\rightarrow D^{r}$ such that
\begin{equation*}
\lambda \circ \xi _{i}=\lambda _{i},\text{ for every }i\in
\mathbb{N.}
\end{equation*}
In particular $\lambda \circ \xi _{r}=\lambda
_{r}=\mathrm{Id}_{D^{r}}.$
\end{proof}

\begin{theorem}
\label{teo: Teo1}Let $(D,\Delta ,\varepsilon )$ be a coalgebra in
a cocomplete abelian monoidal category $\M$ and let $(M,\rho
_{M}^{r},\rho _{M}^{l})$ be a $D$-bicomodule. Let $\delta
:D\rightarrow E$ be a monomorphism which is a morphism of
coalgebras in $\M$. Let $M:=E/D\cot _{E}D\simeq D^2/D.$ Assume
that

i) $D$ is a formally smooth coalgebra in $\M.$

ii) $M$ is $\mathcal{I}$-injective where $\mathcal{I}:=\{f\in {^{D}\M^{D}}%
\mid f$ cosplits in $\M\}.$

iii) For any $r\in \mathbb{N},$ the canonical injection $\xi
_{r}^{r+1}:D^{r}\rightarrow D^{r+1}$ cosplits in $\M.$\newline
Then there is a coalgebra homomorphism
$f_{D}:\widetilde{D}\rightarrow D$ such that
\begin{equation*}
f_{D}\circ \xi _{1}=\mathrm{Id}_{D}
\end{equation*}
and a $D$-bicomodule homomorphism $f_{M}:\widetilde{D}\rightarrow
M$ such that
\begin{equation*}
f_{M}\circ \xi _{2}=\eta (D,D,0).
\end{equation*}
Moreover there is a unique morphism $f:\widetilde{D}\rightarrow
T_{D}^{c}(M)$ such that
\begin{equation}\label{formula: fxin}
f\xi _{n}=\sigma _{n}f_{n},\text{ for any }n\in \N,
\end{equation}
where
\begin{equation*}
f_{n}=\sum_{t=0}^{n}i_{t}^{n}f_{M}^{\square_D t}\overline{\Delta }_{\widetilde{%
D}}^{t-1}\xi _{n}.
\end{equation*}
Also:

1) $f$ is a coalgebra homomorphism;

2) $p_{0}f=f_{D}$ and $p_{1}f=f_{M}$, where
$p_{n}:T_{D}^{c}(M)\rightarrow M^{\square n}$ denotes the
canonical projection. \newline Furthermore, any coalgebra
homomorphism $f:\widetilde{D}\rightarrow T_{D}^{c}(M)$ that
fulfils 2) satisfies the following relation:
\begin{equation}\label{formula: pnf}
p_{k}f=f_{M}^{\square k}\overline{\Delta
}_{\widetilde{D}}^{k-1}\text{ for any }k\in \N.
\end{equation}
\begin{equation*}
\xymatrix@C=2cm{
  T_{D}^{c}(M) \ar[d]_{p_0} \ar[r]^{p_1} & M   \\
  D& \widetilde{D}  \ar[u]_{f_M}\ar@{.>}[ul]|-{f}\ar[l]^{f_D} &D\ar[l]^{\xi_1}  \ar[ul]_0          }
\end{equation*}
\end{theorem}

\begin{proof}
By Theorem \ref{X formally smooth direct limit}, the morphism $\xi
_{1}:D\rightarrow \widetilde{D}$ has a retraction $f_{D}:\widetilde{D}%
\rightarrow D$ which is a coalgebra homomorphism: $f_{D}\circ \xi _{1}=%
\mathrm{Id}_{D}.$ Thus $\widetilde{D}$ becomes a $D$-bicomodule via $%
f_{D}.$ Now we point out that $\xi _{2}\in \mathcal{I}$. In fact
by Lemma \ref{lem: lambda family} $\xi _{2}$ cosplits in $\M.$
Moreover, $D^{2}$ is a
$D$-bicomodule via $f_{D}$ and $\xi _{2}$ a morphism of $D$%
-bicomodules. As explained in (\ref{claim: golden claim}) $(M,\eta (D,D,0))=%
\C(\xi _{1}^{2})$ is a $D$-bicomodule and $\eta (D,D,0)$ is a
morphism of
$D$-bicomodules. Since $\xi _{2}\in \mathcal{I}$ and $M$ is $\mathcal{I}$%
-injective, there exists a $D$-bicomodule homomorphism $f_{M}:\widetilde{%
D}\rightarrow M$ such that
\begin{equation*}
f_{M}\circ \xi _{2}=\eta (D,D,0).
\end{equation*}
Finally $f_{M}\circ \xi _{1}=f_{M}\circ \xi _{2}\circ \xi
_{1}^{2}=\eta (D,D,0)\circ \xi _{1}^{2}=0,$ so that, by Theorem
\ref{teo: pre univ property of cotensor coalgebra}, there exists a
coalgebra homomorphism $f:\widetilde{D}\rightarrow T_{D}^{c}(M)$
which fulfills the required conditions.
\end{proof}

\newpage

\begin{definition}Keep the hypothesis and notations of Theorem \ref{teo:
Teo1}. Then we have:

\begin{equation*}
f_{M}^{\cot _{D}n}\circ \overline{\triangle
}_{\widetilde{D}}^{n-1}\circ \xi _{n}\overset{(\ref{formula:
pnf})}{=}p_{n}\circ f\circ \xi _{n}\overset{(\ref{formula:
fxin})}{=}p_{n}\circ \sigma _{n}\circ f_{n}=p_{n}^{n}\circ
f_{n}=0.
\end{equation*}
Therefore, given any $k\geq 0,$ we have
\begin{equation*}
f_{M}^{\cot _{D}n}\circ \overline{\triangle
}_{\widetilde{D}}^{n-1}\circ \xi
_{n+k}\circ \xi _{n}^{n+k}=f_{M}^{\cot _{D}n}\circ \overline{\triangle }_{%
\widetilde{D}}^{n-1}\circ \xi _{n}=0.
\end{equation*}
Since $[E/D^{n}\cot _{E}D^{k},\eta (D^{n},D^{k},0)]=\C(\xi
_{n}^{n+k}),$ there exists a unique morphism $$\Phi
(D^{n},D^{k},0):E/D^{n}\cot _{E}D^{k}\rightarrow (E/D\cot
_{E}D)^{\cot _{D}n}$$ such that
\begin{equation}
\Phi (D^{n},D^{k},0)\circ \eta (D^{n},D^{k},0)=f_{M}^{\cot
_{D}n}\circ \overline{\triangle }_{\widetilde{D}}^{n-1}\circ \xi
_{n+k}.  \label{def Phi}
\end{equation}
Note that, since $\eta (D^{n},D^{k},0)$, $f_{M}^{\cot _{D}n}\circ \overline{%
\triangle }_{\widetilde{D}}^{n-1}\circ \xi _{n+k}$ and $\xi
_{n}^{n+k}$ are morphisms of $D$-bicomodules, so is $\Phi
(D^{n},D^{k},0).$\newline Moreover observe that
\begin{equation*}
\Phi (D,D,0)\circ \eta (D,D,0)=f_{M}^{\cot _{D}1}\circ \overline{\triangle }%
_{\widetilde{D}}^{0}\circ \xi _{2}=f_{M}\circ \xi _{2}=\eta
(D,D,0).
\end{equation*}
Since $\eta (D,D,0)$ is an epimorphism we get
\begin{equation}\label{formula: Fi trivial}
\Phi (D,D,0)=\mathrm{Id}_{E/D\cot _{E}D}.
\end{equation}
\end{definition}

The proof of the following theorem requires some technicalities
that, for an easier reeding, were included in Appendix
\ref{appendix2}.

\begin{theorem}
\label{teo: super Phi} Let $\mathcal{M}$ be a cocomplete abelian
monoidal category. Keep the hypothesis and notations of Theorem
\ref{teo: Teo1}. Then $\Phi (D^{n},D,0)$ is a monomorphism for any
$n\in\N .$ Moreover, if we assume that the morphism
$E/D^{n}\square _{E}p_{D}^{E}$ is an epimorphism for every $n\in
\N ,$ we get that $\Phi (D^{n},D,0)$ is an isomorphism for any
$n\in \N .$
\end{theorem}

\begin{proof}
By Theorem \ref{teo: induction Phi}, we have (\ref{formula:
induction Phi}):
\begin{equation*}
\Phi (D^{n},D,0)=[\Phi (D^{n-1},D,0)\square _{D}(E/D\square
_{E}D)]\circ \lbrack E/D^{n-1}\square _{E}{^{D}\overline{\rho
}_{E/D\square _{E}D}}]\circ (\gamma_n \square _{E}D).
\end{equation*}%
We point out that, by Proposition \ref{pro: rho bar is iso}, the morphism ${%
^{D}\overline{\rho }_{E/D\square _{E}D}:}E/D\square
_{E}D\rightarrow D\square _{D}(E/D\square _{E}D)$ is always an
isomorphism.\\ As explained in Definition \ref{def: def gamma} the
morphism $\gamma_n :E/D^{n}\rightarrow E/D^{n-1}\square _{E}E/D$
defined by relation
\begin{equation*}
\gamma_n \circ p_{D^{n}}^{E}=(p_{D^{n-1}}^{E}\square
_{E}p_{D}^{E})\circ \overline{\triangle }_{E}.
\end{equation*}%
is always a monomorphism. By Proposition \ref{pro: E/F exact}, $%
(p_{D^{n-1}}^{E}\square _{E}p_{D}^{E})\circ \overline{\triangle
}_{E}$ is an epimorphism whenever $E/D^{n-1}\square _{E}p_{D}^{E}$
is an epimorphism. In this case $\gamma_n $ is an epimorphism too
and hence it is an isomorphism. Since $\Phi (D,D,0)=Id_{E/D\square
_{E}D},$ by induction, we conclude.
\end{proof}

\begin{theorem}
\label{teo: Teo2}Let $(D,\Delta ,\varepsilon )$ be a coalgebra in
a
cocomplete abelian monoidal category $\mathcal{M}$ satisfying AB5. Let $%
\delta :D\rightarrow E$ be a monomorphism which is a morphism of
coalgebras in $\mathcal{M}$. \newline Let $M:=E/D\square
_{E}D\simeq D^2/D.$ Assume that

1) $D$ is a formally smooth coalgebra in $\mathcal{M}.$

2) $M$ is $\mathcal{I}$-injective where $\mathcal{I}:=\{f\in
{^{D}\mathcal{M} ^{D}}\mid f$ cosplits in $\mathcal{M}\}$.

3) For any $r\in \mathbb{N},$ the canonical injection $\xi
_{r}^{r+1}:D^{r}\rightarrow D^{r+1}$ cosplits in
$\mathcal{M}$.\newline Let
$$f:\widetilde{D}\rightarrow T_{D}^{c}(M)$$ be the unique morphism,
constructed in Theorem \ref{teo: Teo1}, such that
\begin{equation*}
f\xi _{n}=\sigma _{n}f_{n},\text{ for any }n\in \mathbb{N},
\end{equation*}
where
\begin{equation*}
f_{n}=\sum_{t=0}^{n}i_{t}^{n}f_{M}^{\square_D t}\overline{\Delta }_{%
\widetilde{D}}^{t-1}\xi _{n}.
\end{equation*}
Then $f$ is a monomorphism in $\mathcal{M}$. Moreover, if we
assume that

4) the morphism $E/D^{n}\square _{E}p_{D}^{E}$ is an epimorphism for every $%
n\in
\mathbb{N}
,$\newline then $f$ is an isomorphism and, with the further
assumptions that

5) $\mathcal{M}$ is complete with right exact tensor functors;

6) denumerable coproducts commute with $\otimes $; \newline then
$\widetilde{D}$ is a formally smooth coalgebra in $\mathcal{M}.$
\end{theorem}

\begin{proof}
Since $\mathcal{M}$ fulfills AB5 condition, it is enough to prove
that $f_{n} $ is a monomorphism (resp. isomorphism) for any $n\geq
0.$ We proceed by induction on $n\geq 0.$

For $n=0$, since $i_{0}^{0}=0$, we have
$f_{n}=f_{0}=i_{0}^{0}f_{M}^{\square
_{D}0}\overline{\Delta }_{\widetilde{D}}^{-1}\xi _{0}=0$. As $%
f_{0}:D^{0}=0\rightarrow C^{0}(M)=0,$ we get that $f_{0}$ is an
isomorphism.

For $n=1$ we have $f_{n}=f_{1}=i_{0}^{1}f_{M}^{\square _{D}0}\overline{%
\Delta }_{\widetilde{D}}^{-1}\xi _{1}+i_{1}^{1}f_{M}^{\square _{D}1}%
\overline{\Delta }_{\widetilde{D}}^{0}\xi _{1}=i_{0}^{1}f_{M}^{\square _{D}0}%
\overline{\Delta }_{\widetilde{D}}^{-1}\xi _{1}=i_{0}^{1}f_{D}\xi
_{1}=i_{0}^{1}\circ \mathrm{Id}_{D}=i_{0}^{1}=\mathrm{Id}_{D},$ as $%
i_{1}^{1}=0.$ Thus $f_{1}$ is an isomorphism.\newline Let $n>1$.
Assume that $f_{n-1}$ is a monomorphism (resp. isomorphism) and
let us prove that $f_{n}$ is a monomorphism (resp. isomorphism)
too. \newline We have
\begin{equation*}
\Phi (D^{n-1},D,0)\circ \eta (D^{n-1},D,0)\overset{(\ref{def Phi})}{=}%
f_{M}^{\square _{D}n-1}\circ \overline{\Delta
}_{\widetilde{D}}^{n-2}\circ
\xi _{n}\overset{(\ref{formula: pnf})}{=}p_{n-1}\circ f\circ \xi _{n}\overset%
{(\ref{formula: fxin})}{=}p_{n-1}\circ \sigma _{n}\circ
f_{n}=p_{n}^{n-1}\circ f_{n},
\end{equation*}%
so that the following diagram
\begin{equation*}
\xymatrix@C=2cm{
  0 \ar[r] & D^{n-1} \ar[d]_{f_{n-1}} \ar[r]^{\xi_{n-1}^n} & D^{n} \ar[d]_{f_n}
  \ar[r]^{\eta(D^{n-1},D,0)} & \frac{E}{D^n}\cot_E D \ar[d]|{\Phi(D^{n-1},D,0)} \ar[r]& 0 \\
  0 \ar[r] & C^{n-1}(M) \ar[r]^{\sigma_{n-1}^n} & C^{n}(M) \ar[r]^{p_n^{n-1}} & M^{\cot_D n-1} \ar[r] & 0   }
\end{equation*}%
commutes and hence $\Phi (D^{n-1},D,0)=\frac{f_{n}}{f_{n-1}}.$
\newline By Theorem \ref{teo: super Phi}, the morphism $\Phi
(D^{n-1},D,0)$ is a
monomorphism (resp. an isomorphism if we assume also that the morphism $%
E/D^{n}\square _{E}p_{D}^{E}$ is an epimorphism for every $n\in
\mathbb{N}
$). Hence, by applying $5$-Lemma to the diagram above, we get that
$f_{n}$ is a monomorphism (resp. isomorphism) too.\newline We
conclude by observing that, if 5) and 6) hold true, by Theorem
\ref{teo Cosmooth}, we get that $T_{D}^{c}(M)$ is formally smooth.
\end{proof}
Let $E$ be a coalgebra in $\mathcal{M}$ and set%
\begin{eqnarray*}
\mathcal{I}^{E} &=&\left\{ f\in \mathcal{M}^{E}\mid f\text{ cosplits in }%
\mathcal{M}\right\} \text{,} \\
^{E}\mathcal{I} &=&\left\{ f\in ^{E}\mathcal{M}\mid f\text{ cosplits in }%
\mathcal{M}\right\} \text{,}
\end{eqnarray*}%
so that%
\begin{equation*}
\mathcal{I=}^{E}\mathcal{I}\cap \mathcal{I}^{E}
\end{equation*}

\begin{lemma}
\label{lem: inj => flat}Let $\mathcal{M}$ be an abelian monoidal
category,
let $N$ be an $\mathcal{I}^{E}$-injective right $E$-comodule and let $%
i_{L}^{M}:L\rightarrow M$ be a morphism in $^{E}\mathcal{I}.$ Then
$N\square _{E}p_{L}^{M}$ is an epimorphism.
\end{lemma}

\begin{proof}
Since the tensor products are left exact, the sequence%
\begin{equation*}
0\rightarrow N\square _{E}L\overset{N\square _{E}i_{L}^{M}}{\rightarrow }%
N\square _{E}M\overset{N\square _{E}p_{L}^{M}}{\rightarrow
}N\square _{E}M/L
\end{equation*}%
is exact. We have to prove that $N\square _{E}p_{L}^{M}$ is an
epimorphism. Since $i_{L}^{M}$ $\in ^{E}\mathcal{I}$, the sequence
\begin{equation*}
0\rightarrow L\overset{i_{L}^{M}}{\rightarrow }M\overset{p_{L}^{M}}{%
\rightarrow }M/L\rightarrow 0
\end{equation*}%
splits in ${\mathcal{M}}$ and hence the sequence%
\begin{equation*}
0\rightarrow N\otimes L\overset{N\otimes i_{L}^{M}}{\rightarrow }N\otimes M%
\overset{N\otimes p_{L}^{M}}{\rightarrow }N\otimes M/L\rightarrow
0
\end{equation*}%
is exact so that the morphism%
\begin{equation*}
\left( N\otimes E\right) \square _{E}p_{L}^{M}:\left( N\otimes
E\right) \square _{E}M\rightarrow \left( N\otimes E\right) \square
_{E}M/L
\end{equation*}%
is an epimorphism since $\left( N\otimes E\right) \square
_{E}p_{L}^{M}\cong
N\otimes p_{L}^{M}$. Now, as $N$ is an $\mathcal{I}^{E}$-injective right $E$%
-comodule, the monomorphism%
\begin{equation*}
\rho _{N}:N\rightarrow N\otimes E
\end{equation*}%
which belongs to $\mathcal{I}^{E}$, has a retraction $\beta _{N}\in {%
\mathcal{M}^{E}}$ and we have%
\begin{equation*}
\left( N\square _{E}p_{L}^{M}\right) \left( \beta _{N}\square
_{E}M\right) =\beta _{N}\square _{E}p_{L}^{M}=\left( \beta
_{N}\square _{E}M/L\right) \left[ \left( N\otimes E\right) \square
_{E}p_{L}^{M}\right] \text{.}
\end{equation*}%
Now $\beta _{N}\square _{E}M/L$ is an epimorphism since it has a section $%
\rho _{N}\square _{E}M/L$, and $\left( N\otimes E\right) \square
_{E}p_{L}^{M}$ is an epimorphism so that we deduce that $N\square
_{E}p_{L}^{M}$ is an epimorphism.
\end{proof}

\begin{lemma}
\label{lem: D alla n}Let $(D,\Delta ,\varepsilon )$ be a coalgebra
in a
cocomplete abelian monoidal category $\mathcal{M}$ satisfying AB5. Let $%
\delta :D\rightarrow E$ be a monomorphism which is a morphism of
coalgebras
in $\mathcal{M}$. Then%
\begin{equation*}
D^{\wedge _{E}n}=D^{\wedge _{\widetilde{D}}n}
\end{equation*}%
for every $n\in
\mathbb{N}
.$
\end{lemma}

\begin{proof}
We proceed by induction on $n\geq 1$ being the cases $n=0,1$
trivial. Assume
that $D^{\wedge _{E}n}=D^{\wedge _{\widetilde{D}}n}.$ Note that%
\begin{equation*}
\ker \left[ \left( p_{D^{\wedge _{E}n}}^{\widetilde{D}}\otimes p_{D}^{%
\widetilde{D}}\right) \Delta _{\widetilde{D}}\right] =D^{\wedge
_{E}n}\wedge
_{\widetilde{D}}D=D^{\wedge _{\widetilde{D}}n}\wedge _{\widetilde{D}%
}D=D^{\wedge _{\widetilde{D}}n+1}.
\end{equation*}%
Let us prove that
\begin{equation*}
\left( D^{\wedge _{E}n+1},\xi _{n+1}\right) =\ker \left[ \left(
p_{D^{\wedge
_{E}n}}^{\widetilde{D}}\otimes p_{D}^{\widetilde{D}}\right) \Delta _{%
\widetilde{D}}\right] .
\end{equation*}%
Let $f:X\rightarrow \widetilde{D}$ be a morphism such that
\begin{equation*}
\left( p_{D^{\wedge _{E}n}}^{\widetilde{D}}\otimes p_{D}^{\widetilde{D}%
}\right) \circ \Delta _{\widetilde{D}}\circ f=0.
\end{equation*}%
We have to prove that there exists a unique morphism $\overline{f}%
:X\rightarrow D^{\wedge _{E}n+1}$ such that%
\begin{equation*}
\xi _{n+1}\circ \overline{f}=f.
\end{equation*}%
We have%
\begin{eqnarray*}
\left( p_{D^{\wedge _{E}n}}^{E}\otimes p_{D}^{E}\right) \circ
\Delta _{E}\circ \widetilde{\delta }\circ f &=&\left( p_{D^{\wedge
_{E}n}}^{E}\otimes p_{D}^{E}\right) \circ \left( \widetilde{\delta
}\otimes
\widetilde{\delta }\right) \circ \Delta _{\widetilde{D}}\circ f \\
&=&\left( \frac{\widetilde{\delta }}{D^{\wedge _{E}n}}\otimes \frac{%
\widetilde{\delta }}{D}\right) \circ \left( p_{D^{\wedge _{E}n}}^{\widetilde{%
D}}\otimes p_{D}^{\widetilde{D}}\right) \circ \Delta
_{\widetilde{D}}\circ f=0.
\end{eqnarray*}%
Since
\begin{equation*}
\left( D^{\wedge _{E}n+1},\delta _{n+1}\right) =\ker \left[ \left(
p_{D^{\wedge _{E}n}}^{E}\otimes p_{D}^{E}\right) \circ \Delta
_{E}\right] ,
\end{equation*}%
by the universal property of the kernel, there exists a unique morphism $%
\overline{f}:X\rightarrow D^{\wedge _{E}n+1}$ such that%
\begin{equation*}
\delta _{n+1}\circ \overline{f}=\widetilde{\delta }\circ f.
\end{equation*}%
Since $\delta _{n+1}=\widetilde{\delta }\circ \xi _{n+1}$ and
since, by AB5, $\widetilde{\delta }$ is a monomorphism, we
conclude.
\end{proof}

\begin{definitions}
Let $C$ be a coalgebra in an abelian monoidal category
$\mathcal{M}$. We say
that a quotient $M/L$ of a right $C$-comodule $M$ is a $\mathcal{I}^{C}$%
\emph{-quotient} of $M,$ whenever the canonical injection $L\rightarrow M$ is in $%
\mathcal{I}^{C}\ $that is the sequence%
\begin{equation*}
0\rightarrow L\rightarrow M\rightarrow \frac{M}{L}\rightarrow 0
\end{equation*}%
is $\mathcal{I}^{C}$-exact.\newline
A coalgebra $C$ is called (right) \emph{hereditary} whenever every $\mathcal{%
I}^{C}$-quotient of an $\mathcal{I}^{C}$-injective comodule is $\mathcal{I}%
^{C}$-injective.
\end{definitions}

\begin{definition}
A coalgebra $C$ is called \emph{coseparable} whenever
$\mathbf{H}^{1}\left( M,C\right) =0$ for every $M\in
{^{C}\mathcal{M}^{C}.}$
\end{definition}

The following theorem was proved for the abelian monoidal category
of vector spaces in \cite{JLMS}.

\begin{theorem}
\label{teo: hered => fs}Let $D$ be a coalgebra in a
cocomplete abelian monoidal category $\mathcal{M}$ satisfying AB5. Let $%
\delta :D\rightarrow E$ be a monomorphism which is a morphism of
coalgebras in $\mathcal{M}$. Assume that

1) $D$ is a coseparable coalgebra in $\mathcal{M}.$

2) For any $r\in \mathbb{N},$ the canonical injection $\xi
_{r}^{r+1}:D^{\wedge _{E}r}\rightarrow D^{\wedge _{E}r+1}$ cosplits in $%
\mathcal{M}$ {(e.g. $\mathcal{M}$ is a semisimple category).}

3) $\widetilde{D}$ is a right hereditary coalgebra in
$\mathcal{M}$.\newline Let
\begin{equation*}
f:\widetilde{D}\rightarrow T_{D}^{c}(\frac{D\wedge _{E}D}{D})
\end{equation*}%
be the unique morphism, constructed in Theorem \ref{teo: Teo1}.
\newline Then $f$ is an isomorphism in $\mathcal{M}$. Moreover, if
we assume that

4) $\mathcal{M}$ is complete with right exact tensor functors;

5) denumerable coproducts commute with $\otimes $; \newline then
$\widetilde{D}$ is a formally smooth coalgebra in $\mathcal{M}.$
\end{theorem}

\begin{proof}
First of all let us point out that, by Lemma \ref{lem: D alla n} one has%
\begin{equation*}
D^{\wedge _{E}n}=D^{\wedge _{\widetilde{D}}n}.
\end{equation*}%
Set $D^{n}=D^{\wedge _{E}n}=D^{\wedge _{\widetilde{D}}n}.$

We apply Theorem \ref{teo: Teo2} to the case
$D=D,E=\widetilde{D},\delta =\xi _{1}.$ In fact, since $D$ is
coseparable, it is in particular a formally smooth coalgebra in
$\mathcal{M}$ \cite[Corollary 4.18]{AMS}. Moreover every
$D$-bicomodule and in particular $D^{2}/D$ is $\mathcal{I}$-injective \cite[%
Theorem 4.4]{AMS}. It remains to prove that the morphism%
\begin{equation*}
\frac{\widetilde{D}}{D^{n}}\square _{\widetilde{D}}p_{D}^{\widetilde{D}}:%
\frac{\widetilde{D}}{D^{n}}\square
_{\widetilde{D}}\widetilde{D}\rightarrow
\frac{\widetilde{D}}{D^{n}}\square
_{\widetilde{D}}\frac{\widetilde{D}}{D}
\end{equation*}%
is an epimorphism. Since, by Lemma \ref{lem: lambda family}, the
canonical injection $\xi_n:D^n\to \widetilde{D}$ cosplits in $\M$
for every $n\in\N$, we get that $\frac{\widetilde{D}}{D^{n}}$ is an $\mathcal{I}^{\widetilde{D}%
}$-quotient of $\widetilde{D}$. Since, by assumption,
$\widetilde{D}$ is a right hereditary, then
$\frac{\widetilde{D}}{D^{n}}$ is $\mathcal{I}^{\widetilde{D}%
}$-injective for every $n\in
\mathbb{N}
.$ Now, $i_{D}^{\widetilde{D}}=\xi _{1}:D\rightarrow
\widetilde{D}$ has a
retraction $f_{D}:\widetilde{D}\rightarrow D$ in $\mathcal{M}$ and hence $%
i_{D}^{\widetilde{D}}\in ^{\widetilde{D}}\mathcal{I}$. By Lemma
\ref{lem: inj => flat}, we conclude.
\end{proof}

\begin{theorem}\label{teo: fs ad hered}
Let $\left( C,\Delta ,\varepsilon \right) $ be a formally smooth
coalgebra in an abelian monoidal category $\mathcal{M}$. Then $C$
is a right hereditary coalgebra.
\end{theorem}

\begin{proof}
Let $M\in {\mathcal{M}}^{C}$ and $\left( \mho _{1}C,\pi \right) =\mathrm{%
Coker}\left( \Delta \right) $. Let us consider the following exact sequence in $%
^{C}{\mathcal{M}}^{C}$:%
\begin{equation*}
0\rightarrow C\overset{\Delta }{\rightarrow }C\otimes C\overset{\pi }{%
\rightarrow }\mho _{1}C\rightarrow 0.
\end{equation*}%
Clearly $r_{C}\left( C\otimes \varepsilon \right) \in
{^{C}\mathcal{M}}$ is
a retraction of $\Delta $ so that the sequence above splits in $^{C}{%
\mathcal{M}}$ and hence
\begin{equation*}
0\rightarrow M\square _{C}C\overset{M\square _{C}\Delta }{\longrightarrow }%
M\square _{C}C\otimes C\overset{M\square _{C}\pi }{\longrightarrow
}M\square _{C}\mho _{1}C\rightarrow 0
\end{equation*}%
is an exact sequence in ${\mathcal{M}}^{C}$. Since $M\square
_{C}C\cong M,$
we get the exact sequence%
\begin{equation*}
\xymatrix@C=1.5cm{
  0 \ar[r] & M \ar[rr]^{\rho^r_M} && M\otimes C \ar[rr]^{(M\cot_C
\pi)[\overline{\rho}_M^C\ot C]} && M\square_C \mho _{1}C \ar[r] &
0 }
\end{equation*}%
in ${\mathcal{M}}^{C}$. Since $C$ is formally smooth, by
\cite[Corollary 4.21 ]{AMS}, $\mho _{1}C$ is
$\mathcal{I}$-injective, so that , by \cite[Theorem 2.3]{Ar Separable}, there is a morphism of $C$-bicomodules $%
j:\mho _{1}C\rightarrow C\otimes X\otimes C$ that cosplits in $^{C}{\mathcal{%
M}}^{C}.$ Thus%
\begin{equation*}
M\square _{C}j:M\square _{C}\mho _{1}C\rightarrow M\square
_{C}C\otimes X\otimes C\cong M\otimes X\otimes C
\end{equation*}%
is a morphism in ${\mathcal{M}}^{C}$ that cosplits in
${\mathcal{M}}^{C}$. Therefore, by \cite[Theorem 2.3]{Ar
Separable}, $M\square _{C}\mho _{1}C$ is
$\mathcal{I}^{C}$-injective and hence the previous sequence is an
$\mathcal{I}^{C}$-injective resolution of length $1.$ Since
this is true for every $M\in {\mathcal{M}}^{C}$, we conclude that $\mathrm{%
\mathrm{Ext}}_{\mathcal{I}^{C}}^{n}(M,N)=0,$ for every $M,N\in {\mathcal{M}}%
^{C}${$,n\geq 2$ i.e. }that $C$ is a right hereditary coalgebra.
\end{proof}

 \begin{theorem}
\label{teo: hered <=> fs}Let $D$ be a coalgebra in a cocomplete
and complete abelian monoidal category $\mathcal{M}$ satisfying
AB5.\newline Let $\delta :D\rightarrow E$ be a monomorphism which
is a morphism of coalgebras in $\mathcal{M}$.\newline Assume that

\begin{enumerate}
\item[1)] $D$ is a coseparable coalgebra in $\mathcal{M}$.

\item[2)] For every $r\in \mathbb{N} $, the canonical injection
$\xi _{r}^{r+1}:D^{\wedge _{E}r}\rightarrow D^{\wedge _{E}r+1}$
cosplits in $\mathcal{M}$ (e.g. $\mathcal{M}$ is a semisimple
category).

\item[3)] $\mathcal{M}$ has right exact tensor functors.

\item[4)] Denumerable coproducts commute with $\otimes $.
\end{enumerate}

Then the following assertions are equivalent.\newline (i)
$\widetilde{D}$ is a formally smooth coalgebra.\newline (ii)
$\widetilde{D}$ is a right hereditary coalgebra.\newline (ii')
$\widetilde{D}$ is a left hereditary coalgebra.\newline (iii)
$\widetilde{D}\simeq T_{D}^{c}(\frac{D\wedge _{E}D}{D})$ as
coalgebras.
\newline
(iv) $\widetilde{D}\simeq T_{D}^{c}(N)$ as coalgebras, for some $D$%
-bicomodule $N$.
\end{theorem}

\begin{proof}
$\left( i\right) \Rightarrow \left( ii\right) $ follows by Theorem
\ref{teo: fs ad hered}.\newline $\left( ii\right) \Rightarrow
\left( iii\right) $ follows by Theorem \ref{teo: hered =>
fs}.\newline $\left( iii\right) \Rightarrow \left( iv\right) $ is
trivial.\newline $\left( iv\right) \Rightarrow \left( i\right) $.
Since $D$ is coseparable, it
is in particular formally smooth as a coalgebra in $\mathcal{M}$ \cite[%
Corollary 4.18]{AMS}. Moreover every $D$-bicomodule and in
particular $N$ is $\mathcal{I}$-injective (see \cite[Theorem
4.4]{AMS}). We conclude by applying Theorem \ref{teo Cosmooth}.
\end{proof}

\begin{corollary}
\cite{JLMS} Let $K$ be a vector space. Let $E$ be a $K$-coalgebra
with
coseparable coradical $D.$ Then the following assertions are equivalent.%
\newline
(i) $E$ is a formally smooth coalgebra.\newline (ii) $E$ is a
right hereditary coalgebra.\newline (ii') $E$ is a left hereditary
coalgebra.\newline (iii) $E\simeq T_{D}^{c}(\frac{D\wedge
_{E}D}{D})$ as coalgebras. \newline (iv) $E\simeq T_{D}^{c}(N)$ as
coalgebras, for some $D$-bicomodule $N$.
\end{corollary}

\begin{proof}
Since $D$ is the coradical of $E$ is well known that
$E=\widetilde{D}$ (see e.g. \cite[Corollary 9.0.4, page 185]{Sw}).
The conclusion follows by Theorem \ref{teo: hered <=> fs} applied
in the case when $\mathcal{M}$ is the category of vector spaces
over $K.$ We point out that, since, in this case, $\mathcal{M}$ is
a semisimple category, the notions of formally smooth and right
hereditary coalgebras reduce to the classical ones.
\end{proof}

\begin{examples} We now provide a number of examples of abelian monoidal
categories for which our results apply. These categories are all
Grothendieck categories and hence cocomplete and complete abelian
categories satisfying AB5. \smallskip

Let $B$ be a bialgebra over a field $K.$\smallskip

$\bullet $ \emph{The category $_{B}{\mathfrak{M}}=(_{B}{\mathfrak{M}}%
,\otimes _{K},K)$, of all left modules over $B$}. The tensor
$V\otimes W$ of two left $B$-modules is an object in
$_{B}{\mathfrak{M}}$ via the diagonal
action; the unit is $K$ regarded as a left $B$-module via $\varepsilon _{B}$%
.\smallskip

$\bullet $ \emph{The category ${_{B}\mathfrak{M}_{B}}=({_{B}\mathfrak{M}_{B}}%
,\otimes _{K},K)$, of all two-sided modules over $B$}. The tensor
$V\otimes W $ of two $B$-bimodules carries, on both sides, the
diagonal action; the unit is $K$ regarded as a $B$-bimodule via
$\varepsilon _{B}$.\smallskip

$\bullet $ \emph{The category ${^{B}\mathfrak{M}}=({^{B}\mathfrak{M}}%
,\otimes _{K},K)$, of all left comodules over $B$}. The tensor product $%
V\otimes W$ of two left $B$-comodules is an object in
$^{B}{\mathfrak{M}}$ via the diagonal coaction; the unit is $K$
regarded as a left $B$-comodule via the map $k\mapsto 1_{B}\otimes
k$.\smallskip

$\bullet $ \emph{The category ${^{B}\mathfrak{M}^{B}}=({^{B}\mathfrak{M}^{B}}%
,\otimes _{K},K)$ of all two-sided comodules over $B$}. The tensor
$V\otimes W$ of two $B$-bicomodules carries, on both sides, the
diagonal coaction; the unit is $K$ regarded as a $B$-bicomodule
via the maps $k\mapsto 1_{B}\otimes k$ and $k\mapsto k\otimes
1_{B}$.\smallskip

$\bullet $ Let $H$ be a Hopf algebra over a field $K$ with bijective antipode.\\ \emph{The category ${_{H}^{H}\mathcal{YD}}=({_{H}^{H}\mathcal{YD}}%
,\otimes _{K},K)$ of left Yetter-Drinfeld modules over $H$}.
Recall that an
object $V$ in ${_{H}^{H}\mathcal{YD}}$ is a left $H$-module and a left $H$%
-comodule satisfying, for any $h\in H,v\in V$, the compatibility
condition:
\begin{equation*}
\sum ({}h_{(1)}v)_{<-1>}{h_{(2)}}\otimes (h_{(1)}v)_{<0>}=\sum
h_{(1)}v_{<-1>}\otimes h_{(2)}v_{<0>}
\end{equation*}%
where $\Delta _{H}\left( h\right) =\sum h_{(1)}\otimes h_{(2)}$
and $\rho \left( v\right) =\sum v_{<-1>}\otimes v_{<0>}$ denote
the comultiplication of $H$ and the left $H$-comodule structure of
$V$ respectively (we used Sweedler notation).

The tensor product $V\otimes W$ of two Yetter-Drinfeld modules is
an object in ${_{H}^{H}\mathcal{YD}}$ via the diagonal action and
the codiagonal coaction; the unit in $_{H}^{H}\mathcal{YD}$ is $K$
regarded as a left $H$-comodule via the map $x\mapsto 1_{H}\otimes
x$ and as a left $H$-module via the counit $\varepsilon
_{H}$.\smallskip

$\bullet $ \emph{The category $_{Q}{\mathfrak{M}}=(_{Q}{\mathfrak{M}}%
,\otimes _{K},K)$, of all left modules over a quasi-bialgebra $Q$
over a field $K$} (see \cite[Definition XV.1.1, page 368]{Ka}).
\end{examples}


\appendix

\section{The Snake Lemma}\label{section: Snake}
We collect here some technical results that are used in the paper.
Some of them might be found in the literature, nevertheless we
decided to include them here for the reader's sake.

First of all, we need to recall some results related to the so
called "Snake Lemma". Thus let $\M$ be an abelian category and
consider in $\M$ the following commutative diagram with exact rows
\begin{equation*}
\xymatrix@R=25pt@C=30pt{
                          && P \ar[d]_{{k_3}'} \ar[r]^{{\alpha_2^3}'} & \K(f_3) \ar[d]_{k_3} \\
  &A_1 \ar[d]_{f_1} \ar[r]^{\alpha_1^2} & A_2 \ar[d]_{f_2} \ar[r]^{\alpha_2^3} & A_3 \ar[d]_{f_3}\ar[r]&0\\
  0\ar[r]&B_1 \ar[d]_{c_1} \ar[r]^{\beta_1^2} & B_2 \ar[r]^{\beta_2^3} & B_3\\
  &\C(f_1)  }
\end{equation*}
where $P=P(\alpha _{2}^{3},k_{3})$ denotes the pullback of
$(\alpha _{2}^{3},k_{3})$. We have
\begin{equation*}
\beta _{2}^{3}f_{2}{{k_3}'}=f_{3}{\alpha _{2}^{3}k_{3}}^{\prime
}=f_{3}k_{3}{\alpha_2^3}'=0.
\end{equation*}
Since $(B_{1},\beta _{1}^{2})=\ker (\beta _{2}^{3}),$ there exists
a unique morphism $\omega :P\rightarrow B_{1}$ such that
\begin{equation}
\beta _{1}^{2}\omega =f_{2}k_{3}^{\prime }.  \label{formula def w}
\end{equation}
Let $\alpha _{1}^{2}=ib$ be the canonical factorization of $\alpha
_{1}^{2}$ as the composition of a monomorphism
$i:\textrm{Im}(\alpha
_{1}^{2})\rightarrow A_{2}$ and an epimorphism $b:A_{1}\rightarrow \textrm{Im}%
(\alpha _{1}^{2})$. As $\alpha _{2}^{3}ib=\alpha _{2}^{3}\alpha
_{1}^{2}=0$ and $b$ is an epimorphism, we get $\alpha
_{2}^{3}\circ i=0=k_{3}\circ 0.$ By the universal property of $P,$
there exists a unique morphism $\xi :\textrm{Im}(\alpha
_{1}^{2})\rightarrow P$ such that
\begin{equation}
k_{3}^{\prime }\xi =i\text{\qquad and\qquad }{\alpha_2^3}'\xi =0.
\label{formula def xi}
\end{equation}
It is straightforward to check that $(\textrm{Im}(\alpha
_{1}^{2}),\xi )=\ker ({\alpha_2^3}').$ We point out that, since
$\alpha _{2}^{3}$ is an epimorphism, also ${\alpha_2^3}'$ is an
epimorphism and hence $(\ker (f_{3}),{\alpha_2^3}')=\C(\xi ).$
\newline Now, we have:
\begin{equation*}
\beta _{1}^{2}\omega \xi b=f_{2}k_{3}^{\prime }\xi
b=f_{2}ib=f_{2}\alpha _{1}^{2}=\beta _{1}^{2}f_{1}.
\end{equation*}
Since $\beta _{1}^{2}$ is a monomorphism, we deduce that
\begin{equation*}
\omega \xi b=f_{1}
\end{equation*}
and hence $c_{1}\omega \xi b=c_{1}f_{1}=0.$ As $b$ is an
epimorphism, we conclude that $c_{1}\omega \xi =0.$ By the
universal property of $\C(\xi ),$
there exists a unique morphism $\overline{\omega }:\ker (f_{3})\rightarrow \C%
(f_{1})$ such that
\begin{equation}
\overline{\omega }{\alpha_2^3}'=c_{1}\omega . \label{formula def w
bar}
\end{equation}
The morphism $\overline{\omega }$ is usually called
\emph{connecting homomorphism.}\\ In fact it is easy to prove the
existence of morphisms $k_1^2,k_2^3,c_1^2,c_2^3$ such that the
following diagram commutes
\begin{equation*}
\xymatrix@R=25pt@C=30pt{
                          &\K(f_1)\ar[d]^{k_1}\ar@{.>}[r]^{k_1^2}& \K(f_2) \ar[d]^{k_2} \ar@{.>}[r]^{k_2^3} & \K(f_3) \ar[d]^{k_3} \\
  &A_1 \ar[d]_{f_1} \ar[r]^{\alpha_1^2} & A_2 \ar[d]^{f_2} \ar[r]^{\alpha_2^3} & A_3 \ar[d]^{f_3}\ar[r]&0\\
  0\ar[r]&B_1 \ar[d]_{c_1} \ar[r]^{\beta_1^2} & B_2\ar[d]^{c_2} \ar[r]^{\beta_2^3} & B_3\ar[d]^{c_3}\\
  &\C(f_1)\ar@{.>}[r]^{c_1^2}& \C(f_2) \ar@{.>}[r]^{c_2^3} & \C(f_3)  }
\end{equation*}
and one has the following well known result.

\begin{theorem}[Snake Lemma]
  The following sequence is exact:
  \begin{equation}
  \xymatrix@C=0.3cm{
    \K(f_1) \ar[rr]^{k_1^2} && \K(f_2) \ar[rr]^{k_2^3} && \K(f_3) \ar[rr]^{\overline{\omega
      }} && \C(f_1) \ar[rr]^{c_1^2} &&
    \C(f_2) \ar[rr]^{c_2^3} && \C(f_3) }
  \end{equation}
\end{theorem}

\begin{proof}
  It is easy to check that the proof of \cite[Lemma 5, page
  206]{Mac} works also in this more general setting (where $\alpha_1^2$ is not assumed to be a monomorphism and $\beta_2^3$ is not necessarily an epimorphism).
\end{proof}

\begin{proposition}\label{pro: A_1=0}
 Let $\M$ be an abelian category and consider in $\M$ the following
commutative diagram with exact rows
\begin{equation*}
\xymatrix@R=25pt@C=30pt{
                          && P \ar[d]_{{k_3}'} \ar[r]^{{\alpha_2^3}'} & \K(f_3) \ar[d]_{k_3} \\
  &0 \ar[d]_{0} \ar[r]^{0} & A_2 \ar[d]_{f_2} \ar[r]^{\alpha_2^3} & A_3 \ar[d]_{f_3}\ar[r]&0\\
  0\ar[r]&B_1 \ar[d]_{\Id_{B_1}} \ar[r]^{\beta_1^2} & B_2 \ar[r]^{\beta_2^3} & B_3\\
  &B_1  }
\end{equation*}
where $P=P(\alpha _{2}^{3},k_{3})$ denotes the pullback of
$(\alpha _{2}^{3},k_{3})$.\newline Then the connecting
homomorphism $\overline{\omega }$ is uniquely defined by the
following relation:
\begin{equation}
  \beta_1^2 \overline{\omega }=f_2 (\alpha_2^3)^{-1}k_3.
\end{equation}
\end{proposition}

\begin{proof}
By (\ref{formula def w bar}) we have $\overline{\omega }{\alpha _{2}^{3}}%
^{\prime }=c_{1}\omega =\omega $. Therefore we have
\begin{equation*}
\beta _{1}^{2}\overline{\omega }{\alpha _{2}^{3}}^{\prime }=\beta
_{1}^{2}\omega \overset{(\ref{formula def
w})}{=}f_{2}k_{3}^{\prime }=f_{2}(\alpha
_{2}^{3})^{-1}\alpha _{2}^{3}k_{3}^{\prime }=f_{2}(\alpha _{2}^{3})^{-1}k_{3}%
{\alpha _{2}^{3}}^{\prime }.
\end{equation*}
Since ${\alpha _{2}^{3}}^{\prime }$ is an epimorphism, we
conclude. Note that, as $\beta _{1}^{2}$ is a monomorphism, we
have that $\overline{\omega } $ is the unique morphism satisfying
the required relation.
\end{proof}

\section{Technical results to prove Theorem \ref{teo: super Phi}}\label{appendix2}

\begin{proposition}
Let $\mathcal{M}$ be a cocomplete abelian monoidal category. Then,
with the hypothesis and notations of Theorem \ref{teo: Teo1}, the
following relations hold true:
\begin{gather}
^{D}\overline{\rho }_{E/D\square _{E}D}\circ \eta
(D,D,0)=[D\square _{D}\eta (D,D,0)]\circ {^{D}\overline{\rho}
}_{D^{2}}=[f_{D}\xi _{2}\square _{D}\eta
(D,D,0)]\circ \overline{\triangle }_{D^{2}}.  \label{rel 1} \\
\lbrack \xi _{1}^{2}f_{D}\xi _{2}\square _{D}\eta (D,D,0)]\circ \overline{%
\triangle }_{D^{2}}=[D^{2}\square _{D}\eta (D,D,0)]\circ
\overline{\triangle
}_{D^{2}}  \label{rel 2} \\
(E/D\square _{E}\xi _{1}^{n})\circ \eta (D,D,0)=\eta
(D,D^{n},0)\circ \xi
_{2}^{n+1}  \label{rel 3} \\
\lbrack \xi _{1}^{2}\square _{D}(E/D\square _{E}D)]\circ
{^{D}\overline{\rho }_{E/D\square _{E}D}}\circ \eta
(D,D,0)=[D^{2}\square _{D}\eta (D,D,0)]\circ
\overline{\triangle }_{D^{2}}  \label{rel 4} \\
\lbrack \xi _{1}^{2}\square _{D}(E/D\square _{E}\xi _{1}^{n})]\circ {^{D}%
\overline{\rho }_{E/D\square _{E}D}}\circ \eta
(D,D,0)=[D^{2}\square _{D}\eta (D,D^{n},0)\circ \xi
_{2}^{n+1}]\circ \overline{\triangle }_{D^{2}}
\label{rel 5} \\
\eta (D^{n-1},D^{2},D)=[E/D^{n-1}\square _{E}\eta (D,D,0)]\circ
\eta
(D^{n-1},D^{2},0)  \label{rel 6} \\
(E/D^{n-1}\square _{E}\xi _{2}^{n+1})\circ \eta
(D^{n-1},D^{2},0)=(p_{D^{n-1}}\circ \delta _{n+1}\square
_{E}D^{n+1})\circ \overline{\triangle }_{D^{n+1}}  \label{rel 7}
\end{gather}
\end{proposition}

\begin{proof}
1) By (\ref{formula: eta is comodule morph}) applied to the case
\begin{equation*}
F=B=D\text{, }E=E\text{ and }\pi =f_{D}\xi _{2}\text{ }
\end{equation*}%
we get%
\begin{equation*}
L=\frac{D\wedge _{E}D}{D}\text{ \quad and \quad }^{D}\rho
_{\frac{D\wedge _{E}D}{D}}\circ p_{D}^{D^{2}}=(f_{D}\otimes
p_{D}^{D^{2}})\circ \Delta _{D\wedge _{E}D}
\end{equation*}

from which, by identifying $(\frac{D^{2}}{D},p_{D}^{D^{2}})$ with $%
(E/D\square _{E}D,\eta (D,D,0)$, we obtain%
\begin{equation*}
^{D}\rho _{E/D\square _{E}D}\circ \eta (D,D,0)=[D\otimes \eta (D,D,0)]\circ {%
^{D}\rho _{D^{2}}}
\end{equation*}%
which means that $\eta (D,D,0):D^{2}\rightarrow E/D\square _{E}D$
is a morphism of left $D$-comodules. Thus we can  apply
Proposition \ref{lem: passo al bar} in the case
\begin{equation*}
W_{1}=D^{2},W_{2}=E/D\square _{E}D,f=\eta (D,D,0),A=F=E=D,\alpha =i_{A}^{E}=%
\mathrm{Id}_{D},
\end{equation*}%
in order to obtain
\begin{equation*}
^{D}\overline{\rho }_{E/D\square _{E}D}\circ \eta
(D,D,0)=[D\square _{D}\eta (D,D,0)]\circ {^{D}\overline{\rho
}}_{D^{2}}
\end{equation*}%
Note that, in view of \ref{pro: rho chiave} applied to the case%
\begin{equation*}
F=D,A=D^{2},B=0,\alpha =\xi _{1}^{2},\pi =f_{D}\xi _{2},
\end{equation*}
we have $^{D}\rho _{D^{2}}=[f_{D}\xi _{2}\otimes D^{2}]\circ
\triangle _{D^{2}}$ so that $^{D}\overline{\rho
}_{D^{2}}=[f_{D}\xi _{2}\square
_{D}D^{2}]\circ \overline{\triangle }_{D^{2}}$ and hence we get (\ref{rel 1}%
). \newline

2) By (\ref{formula: f_D disappears}), applied to the case%
\begin{equation*}
F=B=D\text{, }E=E\text{, }\pi =f_{D}\xi _{2}\text{,}
\end{equation*}%
we get%
\begin{equation*}
\left( D^{2}\otimes p_{D}^{D^{2}}\right) \circ \Delta
_{D^{2}}=(\xi _{1}^{2}f_{D}\xi _{2}\otimes p_{D}^{D^{2}})\circ
\Delta _{D^{2}}
\end{equation*}%
from which, by identifying $(\frac{D^{2}}{D},p_{D}^{D^{2}})$ with $%
(E/D\square _{E}D,\eta (D,D,0)$, we obtain
\begin{equation*}
\lbrack \xi _{1}^{2}f_{D}\xi _{2}\otimes \eta (D,D,0)]\circ
\triangle _{D^{2}}=[D^{2}\otimes \eta (D,D,0)]\circ \triangle
_{D^{2}}
\end{equation*}%
Therefore, since $\eta (D,D,0):D^{2}\rightarrow E/D\square _{E}D$
is a morphism of left $D$-comodules, one has
\begin{eqnarray*}
&&\chi _{D}(D^{2},E/D\square _{E}D)\circ \lbrack \xi
_{1}^{2}f_{D}\xi
_{2}\square _{D}\eta (D,D,0)]\circ \overline{\triangle }_{D^{2}} \\
&=&[\xi _{1}^{2}f_{D}\xi _{2}\otimes \eta (D,D,0)]\circ \chi
_{D}(D^{2},D^{2})\circ \overline{\triangle }_{D^{2}} \\
&=&[\xi _{1}^{2}f_{D}\xi _{2}\otimes \eta (D,D,0)]\circ \triangle _{D^{2}} \\
&=&[D^{2}\otimes \eta (D,D,0)]\circ \triangle _{D^{2}} \\
&=&[D^{2}\otimes \eta (D,D,0)]\circ \chi _{D}(D^{2},D^{2})\circ \overline{%
\triangle }_{D^{2}} \\
&=&\chi _{D}(D^{2},E/D\square _{E}D)\circ \lbrack D^{2}\square
_{D}\eta (D,D,0)]\circ \overline{\triangle }_{D^{2}}.
\end{eqnarray*}%
Since $\chi _{D}(D^{2},E/D\square _{E}D)$ is a monomorphism we obtain (\ref%
{rel 2}). \newline 3) By applying Proposition \ref{pro: naturality
of eta} in the case
\begin{equation*}
E_{1}=E_{2}=E,F_{1}=F_{2}=B_{1}=D,B_{2}=D^{n},A_{1}=A_{2}=0,e=\mathrm{Id}%
_{E},f=\mathrm{Id}_{D},b=\xi _{1}^{n},a=0,
\end{equation*}%
we obtain(\ref{rel 3}). \newline 4) We have
\begin{eqnarray*}
&&[\xi _{1}^{2}\square _{D}(E/D\square _{E}D)]\circ {^{D}\overline{\rho }}%
_{E/D\square _{E}D}\circ \eta (D,D,0) \\
&\overset{(\ref{rel 1})}{=}&[\xi _{1}^{2}\square _{D}(E/D\square
_{E}D)]\circ \lbrack f_{D}\xi _{2}\square _{D}\eta (D,D,0)]\circ \overline{%
\triangle }_{D^{2}} \\
&=&[\xi _{1}^{2}f_{D}\xi _{2}\square _{D}\eta (D,D,0)]\circ \overline{%
\triangle }_{D^{2}} \\
&\overset{(\ref{rel 2})}{=}&[D^{2}\square _{D}\eta (D,D,0)]\circ \overline{%
\triangle }_{D^{2}}.
\end{eqnarray*}%
Hence we get (\ref{rel 4}). \newline 5) We have
\begin{eqnarray*}
&&[\xi _{1}^{2}\square _{D}(E/D\square _{E}\xi _{1}^{n})]\circ {^{D}%
\overline{\rho }}_{E/D\square _{E}D}\circ \eta (D,D,0) \\
&=&[D^{2}\square _{D}(E/D\square _{E}\xi _{1}^{n})]\circ \lbrack
\xi
_{1}^{2}\square _{D}(E/D\square _{E}D)]\circ {^{D}\overline{\rho }}%
_{E/D\square _{E}D}\circ \eta (D,D,0) \\
&\overset{(\ref{rel 4})}{=}&[D^{2}\square _{D}(E/D\square _{E}\xi
_{1}^{n})]\circ \lbrack D^{2}\square _{D}\eta (D,D,0)]\circ \overline{%
\triangle }_{D^{2}} \\
&=&[D^{2}\square _{D}(E/D\square _{E}\xi _{1}^{n})\circ \eta
(D,D,0)]\circ
\overline{\triangle }_{D^{2}} \\
&\overset{(\ref{rel 3})}{=}&[D^{2}\square _{D}\eta
(D,D^{n},0)\circ \xi _{2}^{n+1}]\circ \overline{\triangle
}_{D^{2}}
\end{eqnarray*}%
Hence we obtain (\ref{rel 5}).\newline 6) By applying Proposition
\ref{pro: diagram eta} in the case
\begin{equation*}
E=E,F=D^{n-1},B=D^{2},A=D,
\end{equation*}%
we obtain \ref{rel 6}:
\begin{equation*}
\eta (D^{n-1},D^{2},D)=[E/D^{n-1}\square _{E}\eta (D,D,0)]\circ
\eta (D^{n-1},D^{2},0)
\end{equation*}%
7) By applying (\ref{formula eta ridotta}) in the case
\begin{equation*}
E=E,F=D^{n-1},B=D^{2},A=0,i_{B}^{F\wedge _{E}B}=\xi
_{2}^{n+1},\J{F}{B}=\delta _{n+1},p_{A}^{F\wedge
_{E}B}=\mathrm{Id}_{D^{n+1}},
\end{equation*}%
we obtain \ref{rel 7}.\newline
\end{proof}

\begin{proposition}Let $\M$ be an abelian
monoidal category. Then we have: {\small\begin{equation}
\label{formula varphi}\lbrack (E/D^{n-1}\cot _{E}\xi _{1}^{2})\cot
_{D}(E/D\cot _{E}\xi _{1}^{n})]\circ \varphi
(D^{n-1},D^{2},D)=[\eta (D^{n-1},D^{2},0)\cot _{D}\eta
(D,D^{n},0){]}\circ \overline{\triangle }_{D^{n+1}},
\end{equation}}
where
\begin{equation}\label{def: phi phi}\varphi (D^{n-1},D^{2},D):=(E/D^{n-1}\cot
_{E}{}^{D}\overline{\rho }_{E/D\cot _{E}D})\circ \eta
(D^{n-1},D^{2},D).
\end{equation}
\end{proposition}

\begin{proof}We have:
\begin{eqnarray*}
&&(E/D^{n-1}\cot _{E}\xi _{2}^{n+1}\cot _{D}E/D\cot
_{E}D^{n})\circ \lbrack (E/D^{n-1}\cot _{E}\xi _{1}^{2})\cot
_{D}(E/D\cot _{E}\xi _{1}^{n})]\circ
\varphi (D^{n-1},D^{2},D) \\
&=&(E/D^{n-1}\cot _{E}\xi _{2}^{n+1}\cot _{D}E/D\cot
_{E}D^{n})\circ \lbrack
(E/D^{n-1}\cot _{E}\xi _{1}^{2})\cot _{D}(E/D\cot _{E}\xi _{1}^{n})]\circ  \\
&&\circ (E/D^{n-1}\cot _{E}{}^{D}\overline{\rho }_{E/D\cot
_{E}D})\circ \eta
(D^{n-1},D^{2},D) \\
&\overset{(\ref{rel 6})}{=}&(E/D^{n-1}\cot _{E}\xi _{2}^{n+1}\cot
_{D}E/D\cot _{E}D^{n})\circ \lbrack E/D^{n-1}\cot _{E}(\xi
_{1}^{2}\cot
_{D}E/D\cot _{E}\xi _{1}^{n})^{D}\overline{\rho }_{E/D\cot _{E}D}]\circ  \\
&&\circ \lbrack E/D^{n-1}\cot _{E}\eta (D,D,0)]\circ \eta
(D^{n-1},D^{2},0)
\\
&=&(E/D^{n-1}\cot _{E}\xi _{2}^{n+1}\cot _{D}E/D\cot
_{E}D^{n})\circ \lbrack
E/D^{n-1}\cot _{E}(\xi _{1}^{2}\cot _{D}E/D\cot _{E}\xi _{1}^{n})^{D}%
\overline{\rho }_{E/D\cot _{E}D}\eta (D,D,0)]\circ  \\
&&\circ \eta (D^{n-1},D^{2},0) \\
&\overset{(\ref{rel 5})}{=}&(E/D^{n-1}\cot _{E}\xi _{2}^{n+1}\cot
_{D}E/D\cot _{E}D^{n})\circ \lbrack E/D^{n-1}\cot _{E}[D^{2}\cot
_{D}\eta (D,D^{n},0)\circ \xi _{2}^{n+1}]\circ \overline{\triangle
}_{D^{2}}]\circ
\eta (D^{n-1},D^{2},0) \\
&=&[E/D^{n-1}\cot _{E}[\xi _{2}^{n+1}\cot _{D}\eta
(D,D^{n},0)\circ \xi
_{2}^{n+1}]\circ \overline{\triangle }_{D^{2}}]\circ \eta (D^{n-1},D^{2},0)%
\\&&\text{ since }\xi _{2}^{n+1}\text{ is a coalgebra homomorphism, we get } \\
&=&[E/D^{n-1}\cot _{E}[D^{n+1}\cot _{D}\eta (D,D^{n},0)]\circ \overline{%
\triangle }_{D^{n+1}}\circ \xi _{2}^{n+1}]\circ \eta (D^{n-1},D^{2},0) \\
&=&[E/D^{n-1}\cot _{E}[D^{n+1}\cot _{D}\eta (D,D^{n},0)]\circ \overline{%
\triangle }_{D^{n+1}}]\circ (E/D^{n-1}\cot _{E}\xi
_{2}^{n+1})\circ \eta
(D^{n-1},D^{2},0) \\
&\overset{(\ref{rel 7})}{=}&[E/D^{n-1}\cot _{E}[D^{n+1}\cot
_{D}\eta (D,D^{n},0)]\circ \overline{\triangle }_{D^{n+1}}]\circ
(p_{D^{n-1}}\circ
\delta _{n+1}\cot _{E}D^{n+1})\circ \overline{\triangle }_{D^{n+1}} \\
&=&[p_{D^{n-1}}\circ \delta _{n+1}\cot _{E}[D^{n+1}\cot _{D}\eta
(D,D^{n},0)]]\circ (D^{n+1}\cot _{E}\overline{\triangle
}_{D^{n+1}})\circ
\overline{\triangle }_{D^{n+1}} \\
&=&[p_{D^{n-1}}\circ \delta _{n+1}\cot _{E}D^{n+1}\cot _{D}\eta (D,D^{n},0){]%
}\circ (\overline{\triangle }_{D^{n+1}}\cot _{E}D^{n+1})\circ \overline{%
\triangle }_{D^{n+1}} \\
&=&[(p_{D^{n-1}}\circ \delta _{n+1}\cot _{E}D^{n+1})\circ \overline{%
\triangle }_{D^{n+1}}\cot _{D}\eta (D,D^{n},0){]}\circ \overline{\triangle }%
_{D^{n+1}} \\
&\overset{(\ref{rel 7})}{=}&[(E/D^{n-1}\cot _{E}\xi
_{2}^{n+1})\circ \eta
(D^{n-1},D^{2},0)\cot _{D}\eta (D,D^{n},0){]}\circ \overline{\triangle }%
_{D^{n+1}} \\
&=&(E/D^{n-1}\cot _{E}\xi _{2}^{n+1}\cot _{D}E/D\cot
_{E}D^{n})\circ \lbrack
\eta (D^{n-1},D^{2},0)\cot _{D}\eta (D,D^{n},0){]}\circ \overline{\triangle }%
_{D^{n+1}}.
\end{eqnarray*}
Since $E/D^{n-1}\cot _{E}\xi _{2}^{n+1}\cot _{D}E/D\cot _{E}D^{n}$
is a monomorphism, we get:
\begin{equation*}
\lbrack (E/D^{n-1}\cot _{E}\xi _{1}^{2})\cot _{D}(E/D\cot _{E}\xi
_{1}^{n})]\circ \varphi (D^{n-1},D^{2},D)=[\eta
(D^{n-1},D^{2},0)\cot _{D}\eta (D,D^{n},0){]}\circ
\overline{\triangle }_{D^{n+1}}.
\end{equation*}
\end{proof}

\begin{lemma}
For any $s\geq k\geq 0,$ we have
\begin{equation}
\Phi (D^{n},D^{s},0)\circ (E/D^{n}\cot _{E}\xi _{k}^{s})=\Phi
(D^{n},D^{k},0).  \label{formula Phi - xi}
\end{equation}
\end{lemma}

\begin{proof}
We apply Proposition \ref{pro: naturality of eta} in the case
\begin{equation*}
E_{1}=E_{2}=E,F_{1}=F_{2}=D^{n},B_{1}=D^{k},B_{2}=D^{s},A_{1}=A_{2}=0,e=%
\mathrm{Id}_{E},f=\mathrm{Id}_{D^{n}},b=\xi _{k}^{s},a=0.
\end{equation*}
in order to obtain:
\begin{equation*}
(E/D^{n}\cot _{E}\xi _{k}^{s})\circ \eta (D^{n},D^{k},0)=\eta
(D^{n},D^{s},0)\circ \xi _{n+k}^{n+s}
\end{equation*}
Therefore we get
\begin{eqnarray*}
&&\Phi (D^{n},D^{s},0)\circ (E/D^{n}\cot _{E}\xi _{k}^{s})\circ
\eta
(D^{n},D^{k},0) \\
&=&\Phi (D^{n},D^{s},0)\circ \eta (D^{n},D^{s},0)\circ \xi _{n+k}^{n+s} \\
&=&f_{M}^{\cot _{D}n}\circ \overline{\triangle
}_{\widetilde{D}}^{n-1}\circ
\xi _{n+s}\circ \xi _{n+k}^{n+s} \\
&=&f_{M}^{\cot _{D}n}\circ \overline{\triangle
}_{\widetilde{D}}^{n-1}\circ \xi _{n+k}=\Phi (D^{n},D^{k},0)\circ
\eta (D^{n},D^{k},0)
\end{eqnarray*}
Since $\eta (D^{n},D^{k},0)$ is an epimorphism, we conclude.
\end{proof}

\begin{proposition}
Let $\M$ be a cocomplete abelian monoidal category. Then, with the
hypothesis and notations of Theorem \ref{teo: Teo1}, we have:
\begin{equation}
\lbrack \Phi (D^{n-1},D,0)\cot _{D}(E/D\cot _{E}D)]\circ \varphi
(D^{n-1},D^{2},D)=f_{M}^{\cot _{D}n}\circ \overline{\triangle }_{\widetilde{D%
}}^{n-1}\circ \xi _{n+1}.  \label{formula Phi - f_M}
\end{equation}
\end{proposition}

\begin{proof}

By $(\ref{def Phi})$ we get
\begin{eqnarray}
\Phi (D^{n-1},D^{2},0)\circ \eta (D^{n-1},D^{2},0) &=&f_{M}^{\cot
_{D}n-1}\circ \overline{\triangle }_{\widetilde{D}}^{n-2}\circ \xi
_{n+1}.
\label{def Phi 1} \\
\Phi (D,D^{n},0)\circ \eta (D,D^{n},0) &=&f_{M}^{\cot _{D}1}\circ \overline{%
\triangle }_{\widetilde{D}}^{0}\circ \xi _{1+n}=f_{M}\circ \xi
_{n+1}. \label{def Phi 2}
\end{eqnarray}
By $(\ref{formula Phi - xi})$ we obtain
\begin{eqnarray}
\Phi (D^{n-1},D^{2},0)\circ (E/D^{n-1}\cot _{E}\xi _{1}^{2})
&=&\Phi
(D^{n-1},D,0).  \label{formula Phi - xi 1} \\
\Phi (D,D^{n},0)\circ (E/D\cot _{E}\xi _{1}^{n}) &=&\Phi (D,D,0).
\label{formula Phi - xi 2}
\end{eqnarray}
Finally we compute:
\begin{eqnarray*}
&&f_{M}^{\cot _{D}n}\circ \overline{\triangle
}_{\widetilde{D}}^{n-1}\circ
\xi _{n+1} \\
&=&[f_{M}^{\cot _{D}n-1}\cot _{D}f_{M}]\circ (\overline{\triangle }_{%
\widetilde{D}}^{n-2}\cot _{D}\widetilde{D})\circ \overline{\triangle }_{%
\widetilde{D}}\circ \xi _{n+1} \\
&=&[f_{M}^{\cot _{D}n-1}\circ \overline{\triangle }_{\widetilde{D}%
}^{n-2}\cot _{D}f_{M}]\circ \overline{\triangle
}_{\widetilde{D}}\circ \xi
_{n+1} \\
&=&[f_{M}^{\cot _{D}n-1}\circ \overline{\triangle }_{\widetilde{D}%
}^{n-2}\circ \xi _{n+1}\cot _{D}f_{M}\circ \xi _{n+1}]\circ \overline{%
\triangle }_{D^{n+1}} \\
&\overset{(\ref{def Phi 1}),(\ref{def Phi 2})}{=}&[\Phi
(D^{n-1},D^{2},0)\eta (D^{n-1},D^{2},0)\cot _{D}\Phi
(D,D^{n},0)\eta
(D,D^{n},0)]\circ \overline{\triangle }_{D^{n+1}} \\
&=&[\Phi (D^{n-1},D^{2},0)\cot _{D}\Phi (D,D^{n},0)]\circ \lbrack
\eta
(D^{n-1},D^{2},0)\cot _{D}\eta (D,D^{n},0)]\circ \overline{\triangle }%
_{D^{n+1}} \\
&\overset{(\ref{formula varphi})}{=}&[\Phi (D^{n-1},D^{2},0)\cot
_{D}\Phi (D,D^{n},0)]\circ \lbrack (E/D^{n-1}\cot _{E}\xi
_{1}^{2})\cot _{D}(E/D\cot
_{E}\xi _{1}^{n})]\circ \varphi (D^{n-1},D^{2},D) \\
&=&[\Phi (D^{n-1},D^{2},0)(E/D^{n-1}\cot _{E}\xi _{1}^{2})\cot
_{D}\Phi
(D,D^{n},0)(E/D\cot _{E}\xi _{1}^{n})]\circ \varphi (D^{n-1},D^{2},D) \\
&\overset{(\ref{formula Phi - xi 1}),(\ref{formula Phi - xi
2})}{=}&[\Phi
(D^{n-1},D,0)\cot _{D}\Phi (D,D,0)]\circ \varphi (D^{n-1},D^{2},D) \\
&\overset{(\ref{formula: Fi trivial})}{=}&[\Phi (D^{n-1},D,0)\cot
_{D}(E/D\cot _{E}D)]\circ \varphi (D^{n-1},D^{2},D).
\end{eqnarray*}
\end{proof}

\begin{lemma}
\label{lem: utile 1} Let $\delta :D\rightarrow E$ be a
monomorphism which is a coalgebra homomorphism in an abelian
monoidal category $\M$. Then we have
\begin{equation*}
(D^{m+n},\delta _{m+n})=\ker [(p_{D^{m}}^{E}\cot
_{E}p_{D^{n}}^{E})\circ \overline{\triangle }_{E}].
\end{equation*}
\end{lemma}

\begin{proof}By Remark \ref{rem utile 1}, we have the following exact
sequence:
\begin{equation*}
\xymatrix@C=1.5cm{
  0 \ar[r] & D^m\wedge_E D^n=D^{m+n} \ar[rr]^(.6){\J{D^m}{D^n}=\delta_{m+n}} && E \ar[rr]^(.4){(p _{D^m}^E\cot_E p^E _{D^n}) \circ \overline{\triangle
 } _{E}} && \frac{E}{D^m}\cot_E \frac{E}{D^n}.}
\end{equation*}
Then we conclude.
\end{proof}

\begin{definition}
\label{def: def gamma}By Lemma \ref{lem: utile 1}, we have
\begin{equation*}
(D^{n},\delta _{n})=\ker [(p_{D^{n-1}}^{E}\cot _{E}p_{D}^{E})\circ \overline{%
\triangle }_{E}]
\end{equation*}
so that
\begin{equation*}
(E/D^{n},p_{D^{n}}^{E})=\C(\delta _{n})=\C[\K((p_{D^{n-1}}^{E}\cot
_{E}p_{D}^{E})\circ \overline{\triangle
}_{E})]=Im[(p_{D^{n-1}}^{E}\cot _{E}p_{D}^{E})\circ
\overline{\triangle }_{E}].
\end{equation*}
Thus there exists a unique morphism $\gamma_n :E/D^{n}\rightarrow
E/D^{n-1}\cot _{E}E/D$ such that
\begin{equation}\label{def: gamma}
\gamma_n \circ p_{D^{n}}^{E}=(p_{D^{n-1}}^{E}\cot
_{E}p_{D}^{E})\circ \overline{\triangle }_{E}.
\end{equation}
Obviously $\gamma_n $ is a monomorphism and is a morphism of
$E$-bicomodules.
\end{definition}

\begin{lemma}
Let $\M$ be an abelian monoidal category. We have
\begin{equation}
(\gamma_n \cot _{E}D)\circ \eta (D^{n},D,0)=\eta
(D^{n-1},D^{2},D). \label{formula gamma}
\end{equation}
\end{lemma}

\begin{proof}
We apply (\ref{formula eta ridotta}) in the case
\begin{equation*}
E=E,F=D^{n},B=D,A=0,i_{B}^{F\w B}=\xi _{1}^{n+1},\J{F}{B}=\delta
_{n+1},p_{A}^{F\w B}=\mathrm{Id}_{D^{n+1}}.
\end{equation*}
in order to obtain:
\begin{equation}
(E/D^{n}\cot _{E}\xi _{1}^{n+1})\circ \eta
(D^{n},D,0)=(p_{D^{n}}^{E}\circ \delta _{n+1}\cot
_{E}D^{n+1})\circ \overline{\triangle }_{D^{n+1}}. \label{rel 7 1}
\end{equation}
Therefore we have
\begin{eqnarray*}
&&(E/D^{n-1}\cot _{E}E/D\cot _{E}\xi _{1}^{n+1})\circ (\gamma_n
\cot
_{E}D)\circ \eta (D^{n},D,0) \\
&=&(\gamma_n \cot _{E}D^{n+1})\circ (E/D^{n}\cot _{E}\xi
_{1}^{n+1})\circ \eta
(D^{n},D,0) \\
&\overset{(\ref{rel 7 1})}{=}&(\gamma_n \cot _{E}D^{n+1})\circ
(p_{D^{n}}^{E}\circ \delta _{n+1}\cot
_{E}D^{n+1})\circ \overline{\triangle }_{D^{n+1}} \\
&=&(\gamma_n \circ p_{D^{n}}^{E}\circ \delta _{n+1}\cot
_{E}D^{n+1})\circ
\overline{\triangle }_{D^{n+1}} \\
&\overset{(\ref{def: gamma})}{=}&[(p_{D^{n-1}}^{E}\cot
_{E}p_{D}^{E})\circ \overline{\triangle }_{E}\circ
\delta _{n+1}\cot _{E}D^{n+1}]\circ \overline{\triangle }_{D^{n+1}} \\
&=&[(p_{D^{n-1}}^{E}\circ \delta _{n+1}\cot _{E}p_{D}^{E}\circ
\delta
_{n+1})\circ \overline{\triangle }_{D^{n+1}}\cot _{E}D^{n+1}]\circ \overline{%
\triangle }_{D^{n+1}} \\
&=&(p_{D^{n-1}}^{E}\circ \delta _{n+1}\cot _{E}p_{D}^{E}\circ
\delta _{n+1}\cot _{E}D^{n+1})\circ (\overline{\triangle
}_{D^{n+1}}\cot
_{E}D^{n+1})\circ \overline{\triangle }_{D^{n+1}} \\
&=&(p_{D^{n-1}}^{E}\circ \delta _{n+1}\cot _{E}p_{D}^{E}\circ
\delta
_{n+1}\cot _{E}D^{n+1})\circ (D^{n+1}\cot _{E}\overline{\triangle }%
_{D^{n+1}})\circ \overline{\triangle }_{D^{n+1}} \\
&=&[E/D^{n-1}\cot _{E}(p_{D}^{E}\circ \delta _{n+1}\cot _{E}D^{n+1})%
\overline{\triangle }_{D^{n+1}}]\circ (p_{D^{n-1}}^{E}\circ \delta
_{n+1}\cot _{E}D^{n+1})\circ \overline{\triangle }_{D^{n+1}} \\
&\overset{(\ref{rel 7})}{=}&[E/D^{n-1}\cot _{E}(p_{D}^{E}\circ
\delta _{n+1}\cot _{E}D^{n+1})\overline{\triangle
}_{D^{n+1}}]\circ (E/D^{n-1}\cot
_{E}\xi _{2}^{n+1})\circ \eta (D^{n-1},D^{2},0) \\
&=&[E/D^{n-1}\cot _{E}(p_{D}^{E}\circ \delta _{n+1}\cot _{E}D^{n+1})%
\overline{\triangle }_{D^{n+1}}\xi _{2}^{n+1}]\circ \eta (D^{n-1},D^{2},0) \\
&=&[E/D^{n-1}\cot _{E}(p_{D}^{E}\circ \delta _{n+1}\xi
_{2}^{n+1}\cot
_{E}\xi _{2}^{n+1})\overline{\triangle }_{D^{2}}]\eta (D^{n-1},D^{2},0) \\
&=&(E/D^{n-1}\cot _{E}E/D\cot _{E}\xi _{2}^{n+1})\circ \lbrack
E/D^{n-1}\cot
_{E}(p_{D}^{E}\circ \delta _{2}\cot _{E}D^{2})\overline{\triangle }%
_{D^{2}}]\eta (D^{n-1},D^{2},0) \\
&\overset{(\ref{rel 7 1})}{=}&(E/D^{n-1}\cot _{E}E/D\cot _{E}\xi
_{2}^{n+1})\circ \lbrack E/D^{n-1}\cot _{E}(E/D\cot _{E}\xi
_{1}^{2})\circ
\eta (D,D,0)]\eta (D^{n-1},D^{2},0) \\
&=&(E/D^{n-1}\cot _{E}E/D\cot _{E}\xi _{2}^{n+1}\circ \xi
_{1}^{2})\circ
\lbrack E/D^{n-1}\cot _{E}\eta (D,D,0)]\eta (D^{n-1},D^{2},0) \\
&=&(E/D^{n-1}\cot _{E}E/D\cot _{E}\xi _{1}^{n+1})\circ \lbrack
E/D^{n-1}\cot
_{E}\eta (D,D,0)]\eta (D^{n-1},D^{2},0) \\
&\overset{(\ref{rel 6})}{=}&(E/D^{n-1}\cot _{E}E/D\cot _{E}\xi
_{1}^{n+1})\circ \eta (D^{n-1},D^{2},D).
\end{eqnarray*}
Since $E/D^{n-1}\cot _{E}E/D\cot _{E}\xi _{1}^{n+1}$ is a
monomorphism, we conclude.
\end{proof}

\begin{theorem}\label{teo: induction Phi}
Let $\M$ be a cocomplete abelian monoidal category. Then, with the
hypothesis and notations of Theorem \ref{teo: Teo1}, we have:
\begin{equation} \Phi (D^{n},D,0)=[\Phi (D^{n-1},D,0)\cot
_{D}(E/D\cot _{E}D)]\circ \lbrack E/D^{n-1}\cot
_{E}{^{D}\overline{\rho }_{E/D\cot _{E}D}}]\circ (\gamma_n \cot
_{E}D).\label{formula: induction Phi}
\end{equation}
\end{theorem}

\begin{proof}
We have
\begin{eqnarray*}
&&\lbrack \Phi (D^{n-1},D,0)\cot _{D}(E/D\cot _{E}D)]\circ \lbrack
E/D^{n-1}\cot _{E}{^{D}\overline{\rho }_{E/D\cot _{E}D}}]\circ
(\gamma_n \cot
_{E}D)\circ \eta (D^{n},D,0) \\
&&\overset{(\ref{formula gamma})}{=}[\Phi (D^{n-1},D,0)\cot
_{D}(E/D\cot
_{E}D)]\circ [E/D^{n-1}\cot _{E}{^{D}\overline{\rho }_{E/D\cot _{E}D}}%
]\circ \eta (D^{n-1},D^{2},D) \\
&&\overset{(\ref{def: phi phi})}{=}[\Phi (D^{n-1},D,0)\cot
_{D}(E/D\cot
_{E}D)]\circ \varphi (D^{n-1},D^{2},D) \\
&&\overset{(\ref{formula Phi - f_M})}{=}f_{M}^{\cot _{D}n}\circ \overline{%
\triangle }_{\widetilde{D}}^{n-1}\circ \xi _{n+1} \\
&&\overset{(\ref{def Phi})}{=}\Phi (D^{n},D,0)\circ \eta
(D^{n},D,0).
\end{eqnarray*}
Since $\eta (D^{n},D,0)$ is an epimorphism, we conclude.
\end{proof}

\begin{proposition}
\label{pro: E/F exact}Let $E$ be a coalgebra in an abelian
monoidal category $\mathcal{M}$, let $X$ be a right coideal and
let $Y$ be a left coideal of $E $ in $\mathcal{M}$. Assume that
the morphism $E/X\square _{E}p_{Y}$ is an
epimorphism. Then we have that the following sequence is exact in $\mathcal{M%
}$:
\begin{equation*}
\xymatrix@C=1.5cm{
  0 \ar[r] & X\wedge_E Y \ar[rr]^{\J{X}{Y}} && E \ar[rr]^{(p _{X}\cot_E p
_{Y}) \circ \overline{\triangle } _{E}} && \frac{E}{X} \square_E
\frac{E}{Y} \ar[r] & 0 }
\end{equation*}
\end{proposition}

\begin{proof}
We have:
\begin{equation*}
(E/X\square _{E}p_{Y})\circ \overline{\rho }_{E/X}^{r}\circ
p_{X}=(E/X\square _{E}p_{Y})\circ (p_{X}\square _{E}E)\circ
\overline{\Delta }_{E}=(p_{X}\square _{E}p_{Y})\circ
\overline{\Delta }_{E},
\end{equation*}%
so that $(p_{X}\square _{E}p_{Y})\circ \overline{\Delta }_{E}$ is
an
epimorphism as a composition of epimorphisms. The conclusion follows by (\ref%
{sequence X wedge E}).
\end{proof}


\begin{thebibliography}{100}

\bibitem[AG]{AG} N. Andruskiewitsch, M. Gra\~{n}a, \emph{Braided Hopf algebras
over non-abelian finite groups}. Colloquium on Operator Algebras
and Quantum Groups (Spanish) (VaquerÝas, 1997). Bol. Acad. Nac.
Cienc. (C\'{o}rdoba) 63 (1999), 45--78.

\bibitem[AS]{AS} N. Andruskiewitsch, H-J. Schneider,  \emph{Pointed Hopf
algebras}. New directions in Hopf algebras, 1--68, Math. Sci. Res.
Inst. Publ., 43, Cambridge Univ. Press, Cambridge, 2002.

\bibitem[AMS1]{Cotensor}  A. Ardizzoni, C. Menini and D. \c{S}tefan, \emph{%
Cotensor Coalgebras in Monoidal Categories}, Comm. Algebra, to
appear.

\bibitem[AMS2]{AMS}  A. Ardizzoni, C. Menini and D. \c{S}tefan, \emph{%
Hochschild Cohomology and 'Smoothness' in Monoidal Categories}, J.
Pure Appl. Algebra, to appear.

\bibitem[Ar]{Ar Separable} A. Ardizzoni,
\textit{Separable Functors and formal smoothness}, submitted.

\bibitem[Ch]{Ch} W. Chin, \emph{Hereditary and path coalgebras}.  Comm. Algebra
\textbf{30} (2002),  no. 4, 1829-1831.

\bibitem[JLMS]{JLMS} P. Jara, D. Llena, L. Merino and D. \c{S}tefan, \emph{Hereditary
and formally smooth coalgebras}, Algebr. Represent. Theory
\textbf{8} (2005), 363-374.

\bibitem[Ka]{Ka}  C. Kassel, \emph{Quantum groups}, Graduate Text in
Mathematics \textbf{155}, Springer, 1995.

\bibitem[Mac]{Mac} S. Mac Lane, \emph{Categories for the working mathematician.
Second edition.} Graduate Texts in Mathematics, 5.
Springer-Verlag, New York, 1998.

\bibitem[Maj]{Maj2}  S. Majid, \emph{Foundations of quantum group theory},
Cambridge University Press, 1995.

\bibitem[Mo]{Mo}  S. Montgomery, \emph{Hopf Algebras and their actions on
rings, }CMBS Regional Conference Series in Mathematics
\textbf{82}, 1993.

\bibitem[Ni]{Ni} Nichols, Warren D. Bialgebras of type one. Comm. Algebra 6
(1978), no. 15, 1521--1552.

\bibitem[Sw]{Sw} M. Sweedler, \emph{Hopf Algebras}, Benjamin, New York, 1969.
\end{thebibliography}
\end{document}